\documentclass[11pt,a4paper]{article}

\usepackage[margin=1in]{geometry}
\usepackage[T1]{fontenc}
\usepackage{lmodern}
\usepackage{amsmath,amssymb,amsthm,mathtools}
\usepackage{microtype}
\usepackage{enumitem}
\usepackage{tikz}
\usepackage{placeins}
\usetikzlibrary{arrows.meta,calc,patterns,positioning,decorations.pathmorphing,decorations.pathreplacing}
\usepackage{xcolor}
\definecolor{linkblue}{RGB}{0,0,128}
\usepackage[colorlinks=true,linkcolor=linkblue,citecolor=linkblue,urlcolor=linkblue]{hyperref}
\usepackage[nameinlink,noabbrev]{cleveref}

\allowdisplaybreaks
\numberwithin{equation}{section}
\setlist{leftmargin=2em,itemsep=0.25em,topsep=0.35em}

\newtheorem{theorem}{Theorem}[section]
\newtheorem{proposition}[theorem]{Proposition}
\newtheorem{lemma}[theorem]{Lemma}
\newtheorem{corollary}[theorem]{Corollary}
\theoremstyle{definition}

\theoremstyle{remark}

\newcommand{\R}{\mathbb R}
\newcommand{\M}{\mathcal M}
\newcommand{\U}{\mathcal U}
\newcommand{\dd}{\mathop{}\!\mathrm d}
\newcommand{\tr}{\operatorname{tr}}
\newcommand{\Ric}{\operatorname{Ric}}
\newcommand{\Scal}{\operatorname{Scal}}
\newcommand{\Rtrap}{R_{\mathrm{trap}}}
\newcommand{\Rcoll}{R_{\mathrm{coll}}}
\newcommand{\Rglue}{R_{\mathrm{glue}}}

\newcommand{\HM}{m}
\newcommand{\msp}{\widetilde m}

\title{A Construction of Cauchy Data for Trapped Surface Formation in the Massive Einstein--Klein--Gordon System}
\author{Xueyi Yu}
\date{November 2026}
\hypersetup{
  pdftitle={A Construction of Cauchy Data for Trapped Surface Formation in the Massive Einstein-Klein-Gordon System},
  pdfauthor={Xueyi Yu}
}

\begin{document}
\maketitle
\vspace{-2\baselineskip}

\begin{abstract}
We construct a class of smooth, spherically symmetric, asymptotically flat
Cauchy data for the massive Einstein--Klein--Gordon system that are initially untrapped but give rise to trapped surface formation in their
maximal Cauchy developments. These data coincide with contracting
Friedmann--Lema\^{i}tre--Robertson--Walker (FLRW) data on an interior ball
and with Schwarzschild data outside a larger ball. The data can also be made to have arbitrarily small compactness: the ratio $2\HM/r$, where $\HM$
is the Hawking mass and $r$ is the area radius, can be made arbitrarily
small throughout the initial Cauchy hypersurface.
\end{abstract}

\section{Introduction}\label{sec:introduction}

Let $(\M,g)$ be a Lorentzian spacetime with signature $(-,+,+,+)$, and
let $\phi$ be a real scalar field on $\M$. The massive Einstein--Klein--Gordon system is
\begin{align}
 \operatorname{Ric}_{\alpha\beta}-\frac12\operatorname{Scal}_g g_{\alpha\beta}
 &=2\left(\partial_\alpha\phi\,\partial_\beta\phi
 -\frac12g_{\alpha\beta}
 \bigl(\partial^\gamma\phi\,\partial_\gamma\phi
       +\mu^2\phi^2\bigr)\right),
 \label{eq:ekg-einstein}\\*
 \Box_g\phi-\mu^2\phi&=0,
 \label{eq:ekg-kg}
\end{align}
where $\mu>0$ is the mass of the scalar field.

A \emph{trapped surface} is a smooth closed spacelike two-surface whose two
future null expansions are everywhere negative. A \emph{weakly outer
trapped surface} is a smooth closed two-sided spacelike two-surface whose
outward future null expansion is everywhere nonpositive. An initial
Cauchy hypersurface is called \emph{initially untrapped} if it contains
no weakly outer trapped surface.

\subsection{Main results}\label{sec:main-results}

We construct smooth, spherically symmetric, asymptotically flat initial
data for the massive Einstein--Klein--Gordon system that are initially untrapped but give rise to the formation of trapped surfaces in their future Cauchy developments. These data coincide
with the data induced on the initial slice of a contracting Friedmann--Lema\^{i}tre--Robertson--Walker
(FLRW) solution on an interior ball. The FLRW solution consists of a
spacetime metric and a spatially homogeneous scalar field:
\begin{equation}\label{eq:flrw-intro}
 g_{\mathrm{FLRW}}=-\dd t^2+a(t)^2
 \bigl(\dd s^2+S_\kappa(s)^2\dd\omega^2\bigr),
 \qquad \phi=\phi(t).
\end{equation}
Here $t$ is comoving proper time, $a(t)$ is the scale factor, normalized
to $a(0)=1$, and $\dd\omega^2$ is the metric of the unit sphere in
angular coordinates. $\kappa$ is the initial spatial curvature, and $s$
is the initial radial geodesic radius. $S_\kappa$ is defined by
\begin{equation}\label{eq:S-kappa-intro}
 S_\kappa(s)=
 \begin{cases}
 \kappa^{-1/2}\sin(\sqrt\kappa\,s),&\kappa>0,\\
 s,&\kappa=0,\\
 (-\kappa)^{-1/2}\sinh(\sqrt{-\kappa}\,s),&\kappa<0.
 \end{cases}
\end{equation}
The initial area radius is $r:=S_\kappa(s)$, since the centered sphere
at $(0,s)$ has area $4\pi r^2$. For $\kappa>0$, $S_\kappa$ and its
inverse are restricted to the increasing branch
$0\le s<\pi/(2\sqrt\kappa)$.

With $H:=\dot a/a$, \eqref{eq:ekg-einstein}--\eqref{eq:ekg-kg} give the
Friedmann and Klein--Gordon equations:
\begin{align}
 H^2+\frac\kappa{a^2}
 &=\frac{\dot\phi^2+\mu^2\phi^2}{3},
 \label{eq:flrw-constraint-intro}\\
 \frac{\ddot a}{a}
 &=\frac{\mu^2\phi^2-2\dot\phi^2}{3},
 \label{eq:flrw-evolution-intro}\\
 \ddot\phi+3H\dot\phi+\mu^2\phi&=0.
 \label{eq:flrw-kg-intro}
\end{align}
Prescribe $\phi(0)=\Phi_0$ and $\dot\phi(0)=\Pi_0$, with
$(\Phi_0,\Pi_0)\ne(0,0)$ and
$(\Pi_0^2+\mu^2\Phi_0^2)/3>\kappa$. Set $H_0:=\dot a(0)$. Recall that $a(0)=1$. From \eqref{eq:flrw-constraint-intro}, we have
\begin{equation}\label{eq:E0-H0}
 H_0=-\sqrt{\frac{\Pi_0^2+\mu^2\Phi_0^2}{3}-\kappa},
\end{equation}
where $H_0$ is chosen on the contracting branch. For fixed mass $\mu>0$
and curvature parameter $\kappa$, these scalar initial data determine
a unique contracting FLRW solution of
\eqref{eq:flrw-constraint-intro}--\eqref{eq:flrw-kg-intro}.
Let $[0,T^*)$ be the maximal interval starting at $t=0$ on which
the solution remains regular and satisfies $a>0$ and $\dot a<0$, and call it the \emph{maximal regular contracting interval}.

The contracting FLRW solution already contains trapped spheres on its
initial slice. To obtain untrapped initial data, we restrict the FLRW
slice to an interior ball $0\le r<r_0$, called the \emph{FLRW ball}.
We choose the ball small enough that it contains no trapped surfaces
and admits an untrapped extension. The extension coincides with a Schwarzschild
slice outside the sphere at $r=R_0$. More precisely, the initial
data $(h,K,\Phi,\Pi)$ satisfy
\begin{equation}\label{eq:completed-data-intro}
 (h,K,\Phi,\Pi)=
 \begin{cases}
 \displaystyle
 \left(\frac{\dd r^2}{1-\kappa r^2}+r^2\dd\omega^2,
 H_0h,\Phi_0,\Pi_0\right),&0\le r<r_0,\\[6pt]
 \displaystyle
 \left(\frac{\dd r^2}{1-2M/r}+r^2\dd\omega^2,0,0,0\right),
 &r>R_0.
 \end{cases}
\end{equation}
Here $h$ is the induced metric,
$K=\frac12\mathcal L_nh$, $\Phi=\phi|_{t=0}$, and
$\Pi=n(\phi)|_{t=0}$, where $n$ is the future unit normal to the initial
slice. The data in the transition region $r_0\le r\le R_0$ satisfy the
constraint equations.

Three threshold area radii appear in our construction: the gluing
threshold $\Rglue$, the trapping threshold $\Rtrap$, and the
collapse threshold $\Rcoll$. These radii are determined by
the corresponding FLRW solution:
\begin{align}
 \Rglue&:=\min\left\{
 \sqrt{\frac{3}{\Pi_0^2+\mu^2\Phi_0^2}},\frac1{\mu|\Phi_0|}
 \right\},
 \label{eq:Rglue-intro}\\
 \Rtrap&:=
 S_\kappa\left(
 \inf_{0<t<T^*}
 \left[
 \int_0^t\frac{\dd\tau}{a(\tau)}
 +S_\kappa^{-1}\left(\frac1{\sqrt{\dot a(t)^2+\kappa}}\right)
 \right]\right).
 \label{eq:thresholds-intro}
\end{align}
The second entry in \eqref{eq:Rglue-intro} is understood as $+\infty$
when $\Phi_0=0$. For $\kappa\le0$, \eqref{eq:flrw-constraint-intro} and
$H(0)<0$ imply that $H(t)<0$ throughout the regular evolution,
so the FLRW solution remains contracting. Thus, if a curvature
singularity forms in finite time, it occurs at the endpoint
$t=T^*$ of the maximal regular contracting interval.
If $T^*<\infty$ and $a^{-1}\in L^1(0,T^*)$, we define
\begin{equation}\label{eq:Rcoll-intro}
 \Rcoll:=
 S_\kappa\left(\int_0^{T^*}\frac{\dd t}{a(t)}\right).
\end{equation}
The importance of these three threshold radii is as follows. An FLRW ball with $r_0<\Rglue$ admits an untrapped completion. The
condition $r_0>\Rtrap$ ensures trapped surface formation
within the future domain of dependence of the FLRW ball. Under the
hypotheses of Corollary~\ref{cor:causal-collapse}, $r_0>\Rcoll$ ensures
curvature blow-up along the central geodesic in this region. Note that
$\Rglue$ depends on $(\mu,\Phi_0,\Pi_0)$, while $\Rtrap$
and $\Rcoll$ depend on $(\mu,\Phi_0,\Pi_0,\kappa)$.
Comparing $\Rglue$ with $\Rtrap$ and $\Rcoll$ gives the following
theorems.

\begin{theorem}\label{thm:general}
Fix $\mu>0$, $\kappa\in\R$, and $(\Phi_0,\Pi_0)$ such that
$(\Pi_0^2+\mu^2\Phi_0^2)/3>\kappa$. Assume that
\begin{equation}\label{eq:kinetic-cone-intro}
 \Pi_0^2>2\mu^2\Phi_0^2.
\end{equation}
Then $\Rtrap<\Rglue$. For every FLRW ball with area radius $r_0$
satisfying
\[
 \Rtrap<r_0<\Rglue,
\]
there exist constants $M>0$ and $R_0>r_0$, and smooth, spherically
symmetric, asymptotically flat Cauchy data
$(\R^3,h,K,\Phi,\Pi)$ satisfying \eqref{eq:completed-data-intro}.
The initial Cauchy hypersurface contains no closed weakly outer trapped
surface, while its maximal Cauchy development contains a trapped
surface and is therefore future null geodesically incomplete.
\end{theorem}

While Theorem~\ref{thm:general} applies to every admissible curvature parameter
$\kappa$, the next theorem considers
the limit $\kappa\to-\infty$.

\begin{theorem}\label{thm:small-radius}
Fix $\mu>0$, $\Phi_0\in\R$, and $\Pi_0\ne0$.
For sufficiently negative $\kappa$, we have $T^*<\infty$ and
$a^{-1}\in L^1(0,T^*)$, so $\Rcoll$ is well defined. We have
\begin{equation}\label{eq:threshold-asymptotic-intro}
 \Rcoll=\Rtrap\longrightarrow0
 \qquad(\kappa\to-\infty),
\end{equation}
whereas $\Rglue$ is independent of $\kappa$. Thus, for all sufficiently
negative $\kappa$, we have
$\Rcoll=\Rtrap<\Rglue$.
For every FLRW ball with area radius $r_0$ satisfying
\[
 \Rcoll=\Rtrap<r_0<\Rglue,
\]
there exist constants $M>0$ and $R_0>r_0$, and smooth, spherically
symmetric, asymptotically flat Cauchy data
$(\R^3,h,K,\Phi,\Pi)$ satisfying \eqref{eq:completed-data-intro}.
The data are initially untrapped, while trapped surfaces form in the
future development. Moreover, along the central comoving geodesic,
$\operatorname{Scal}_g\to-\infty$ in finite proper time, so the Cauchy development is also timelike geodesically incomplete.
\end{theorem}

The next corollary combines Theorem~\ref{thm:small-radius} with the mass-ratio
estimate of Proposition~\ref{prop:small-initial-mass-ratio}.

\begin{corollary}[arbitrarily small compactness]
\label{cor:small-initial-compactness}
Fix $\mu>0$, $r_0>0$, and $\delta\in(0,1/2)$.
There exist initial data as in Theorem~\ref{thm:small-radius}, with an FLRW
ball of area radius $r_0$, such that
\begin{equation}\label{eq:small-compactness-corollary}
 \sup_\Sigma\frac{2\HM}{r}<\delta,
\end{equation}
where $\Sigma$ is the initial hypersurface, $\HM$ is the Hawking mass,
and $r$ is the area radius. The maximal Cauchy development of these
data contains trapped surfaces.
\end{corollary}

The initial compactness $2\HM/r$ can thus be made arbitrarily small
throughout the initial slice compared with the trapping threshold $2\HM/r=1$,
yet trapped surfaces form dynamically in our spacetime.

\begin{figure}[!htbp]
\centering
\begin{tikzpicture}[x=1.45cm,y=1.45cm,>=Latex,
  every node/.style={font=\small},line cap=round,line join=round]
  \coordinate (O) at (0,0);
  \coordinate (B) at (4.4,0);
  \coordinate (C) at (5.8,0);
  \coordinate (I) at (8.6,0);
  \coordinate (S) at (0,3.4);
  \coordinate (Q) at (1,3.4);
  \coordinate (A) at (3.4,1);
  \fill[black!3] (O) -- (I) -- (4.8,3.8)
    .. controls (3.4,3.95) and (2.0,3.4) .. (Q) -- (S) -- cycle;
  \fill[black!6] (O) -- (B) -- (Q) -- (S) -- cycle;
  \fill[black!18] (A) -- (Q) -- (S)
    .. controls (1.0,2.5) and (2.4,1.4) .. (A) -- cycle;
  \draw[black!55,densely dotted,line width=1.05pt]
    (Q) .. controls (2.0,3.4) and (3.4,3.95) .. (4.8,3.8) -- (I);
  \draw[thick] (O) -- (S);
  \draw[line width=1.15pt] (O) -- (I);
  \draw[thick] (B) -- (Q);
  \draw[dashed,thick] (A)
    .. controls (2.4,1.4) and (1.0,2.5) .. (S);
  \draw[line width=1.2pt,decorate,
    decoration={zigzag,segment length=4pt,amplitude=2.3pt}]
    (0.02,3.4) -- (0.98,3.4);
  \node[anchor=south west,inner sep=1pt] at (-0.12,3.51) {curvature blow-up};
  \node[rotate=90,anchor=south] at (-0.10,1.55) {$r=0$};
  \node[align=center] at (1.35,0.63) {$D^+(\U)$};
  \node[rotate=-45,anchor=south] at (3.48,1.01) {$\partial D^+(\U)$};
  \node[align=center,font=\footnotesize] at (0.95,1.82)
    {apparent horizon};
  \draw[black!65,thin] (1.97,1.82) -- (1.69,2.02);
  \filldraw[fill=white,draw=black,line width=1pt] (2.20,2.0) circle (2.8pt);
  \draw[thin] (2.27,2.04) -- (3.28,2.55);
  \node[anchor=west,font=\footnotesize] at (3.28,2.55) {trapped sphere};
  \node[align=center] at (5.47,1.57)
    {Future outside $D^+(\U)$};
  \node[align=center,font=\footnotesize,text=black!65] at (4.20,3.25)
    {undetermined\\future boundary};
  \draw[thin] (B) +(0,-.07) -- +(0,.07);
  \draw[thin] (C) +(0,-.07) -- +(0,.07);
  \fill[pattern=north east lines,pattern color=black!40]
    (4.4,-.045) rectangle (5.8,.045);
  \node[anchor=south west,font=\footnotesize] at (4.20,.08) {$r=r_0$};
  \node[anchor=south west,font=\footnotesize] at (5.40,.08) {$r=R_0$};
  \fill (I) circle (1.5pt);
  \node[anchor=south west] at (8.62,.06) {$i^0$};
  \node[anchor=south] at (7.12,.08) {$\Sigma$};
  \draw[decorate,decoration={brace,mirror,amplitude=4pt,raise=3pt}]
    (0,0) -- (4.4,0) node[midway,below=11pt] {interior FLRW ball $\U$};
  \draw[decorate,decoration={brace,mirror,amplitude=4pt,raise=3pt}]
    (4.4,0) -- (5.8,0)
    node[midway,below=11pt,align=center,font=\footnotesize] {transition\\region};
  \draw[decorate,decoration={brace,mirror,amplitude=4pt,raise=3pt}]
    (5.8,0) -- (8.6,0)
    node[midway,below=11pt,align=center,font=\footnotesize]
    {exterior Schwarzschild region};
\end{tikzpicture}
\caption{Schematic of the spacetime in Theorem~\ref{thm:small-radius}.
Each point away from the center represents a symmetry sphere.
The light gray region is $D^+(\U)$. The dashed curve marks the apparent
horizon, and the darker region contains trapped spheres.
The zigzag segment indicates curvature blow-up. The initial hypersurface
$\Sigma$ is untrapped, and $i^0$ denotes spatial infinity.}
\label{fig:causal-schematic}
\end{figure}
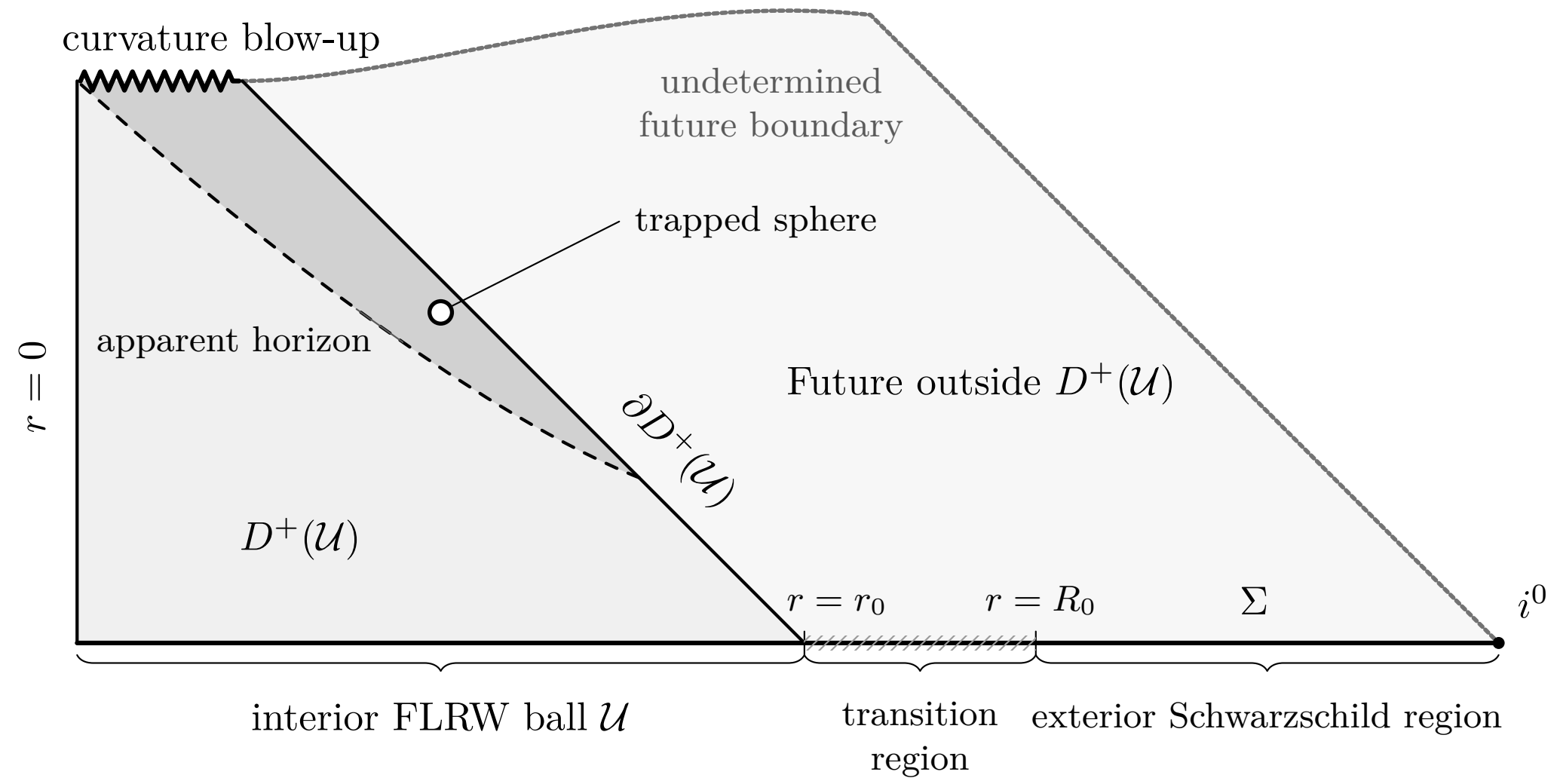

The construction in Theorem~\ref{thm:small-radius}, also used in
Corollary~\ref{cor:small-initial-compactness}, is illustrated in
\cref{fig:causal-schematic}. Trapped spheres and curvature blow-up
occur in the FLRW region determined by the initial ball, while the
full evolution outside this region is not described by our argument.

\FloatBarrier

\subsection{Background and motivation}\label{sec:background}

Trapped surfaces have been used as a sign of gravitational collapse
since Penrose proved his singularity theorem \cite{Penrose1965}. Under
the null convergence condition and with a noncompact Cauchy hypersurface,
the theorem implies that the formation of a trapped surface leads to
future null geodesic incompleteness. We are interested in finding
initially untrapped Cauchy data whose developments contain trapped
surfaces. The answer depends on the matter model.

Christodoulou initiated the study of the evolution of a massless real
scalar field in spherical symmetry
\cite{Christodoulou1986a,Christodoulou1986b,Christodoulou1993}.
In his setting, the initial data are given on an outgoing null
hypersurface. He established sufficient conditions on the initial data for trapped surfaces to form in the future development~\cite{Christodoulou1991}.
Building on Christodoulou's work, Dafermos~\cite{Dafermos2009} constructed
spherically symmetric, asymptotically flat solutions of the
Einstein--massless scalar field system that are causally geodesically
complete to the past and contain black holes in the future.

By Birkhoff's theorem, a spherically symmetric vacuum spacetime is
locally isometric to a Schwarzschild spacetime. To obtain dynamical
trapped surface formation from regular vacuum data, one must go beyond
spherical symmetry.
Christodoulou \cite{Christodoulou2009} achieved this by introducing the
short-pulse method. This construction relies on a double-null
foliation of the spacetime, and the initial data are given on an
outgoing null hypersurface. There is an open set of initial data that
contain no trapped surfaces but lead to their formation. Subsequently,
Klainerman--Rodnianski \cite{KlainermanRodnianski2012} generalized Christodoulou's work by enlarging the admissible set of initial
conditions and reducing the number of derivatives required in the
estimates. Klainerman--Luk--Rodnianski
\cite{KlainermanLukRodnianski2014} extended the results to the case of fully anisotropic
initial data. Li--Yu~\cite{LiYu2015} combined Christodoulou's short-pulse
construction with the gluing method and established trapped surface
formation from smooth, asymptotically flat vacuum Cauchy data. An--Luk~\cite{AnLuk2017} established a
scale-critical criterion for trapped surface formation in vacuum,
allowing weaker incoming radiation than in Christodoulou's construction.

Trapped surface formation has also been studied in other settings,
including charged scalar fields~\cite{AnLim2022} and scalar fields with
a positive cosmological constant~\cite{Costa2020} in spherical symmetry.
Beyond spherical symmetry, examples include the massless Einstein--scalar field~\cite{ZhaoHilditchValienteKroon2024} and Einstein--Maxwell--charged scalar field~\cite{ShenWan2025} systems.

The results listed above concern trapped surface formation for massless
scalar fields or in vacuum. Double-null coordinates are a natural
choice in these settings, as they are adapted to the evolution of the null expansions. The present paper considers
the case of a massive scalar field. The mass term $\mu^2\phi$
in the Klein--Gordon equation \eqref{eq:ekg-kg} introduces an additional
term in the evolution equations in double-null coordinates. Controlling
this term requires a bound for $\phi$ itself as well as its derivatives. This makes it more difficult to close
the estimates used in the massless case.

The rescaling
\begin{equation}\label{eq:covariance}
 (g,\phi,\mu)\longmapsto
 (\lambda^2g,\phi,\mu/\lambda),
 \qquad \lambda>0,
\end{equation}
preserves the massive Einstein--Klein--Gordon system \eqref{eq:ekg-einstein}--\eqref{eq:ekg-kg}.
Since a length scale $L$ transforms to $\lambda L$ and the mass parameter $\mu$ transforms to $\mu/\lambda$ under the rescaling, the dimensionless
product $\mu L$ is invariant. The mass term introduces a length
scale $\mu^{-1}$, which is called the Compton length~\cite{BradyChambersGoncalves1997}. Unlike the
massless equations, the equations at fixed $\mu>0$ are not invariant under
this rescaling. Let $L$ be a characteristic length scale of the spacetime; then the
dimensionless product $\mu L$ matters for the study of collapse.
For example, numerical studies show that the ratio of the initial
pulse width to the Compton length can qualitatively affect the collapse
dynamics~\cite{BradyChambersGoncalves1997,OkawaCardosoPani2014}.

There are existing results for trapped surface and black hole formation
in massive scalar field systems. Giorgi--Shen--Wan \cite{GiorgiShenWan2025} used the
short-pulse method to establish trapped surface formation in the
Einstein--Maxwell--Klein--Gordon system. The scalar field in their
construction can have a quadratic potential
$V(\phi)=\mu^2|\phi|^2$, with $\mu$ sufficiently small so that the potential can be treated as a small
perturbation. This condition can be interpreted in a
scaling-invariant manner as $\mu L\ll1$, where $L$ is the local length
scale normalized to one in their construction. Gon\c{c}alves--Moss
\cite{GoncalvesMoss1997} considered the case $\mu L\gg1$, where $L$
is the characteristic radius of the initial scalar-field distribution.
They discussed black hole formation for a spherically symmetric massive scalar field
using the WKB (Wentzel--Kramers--Brillouin) approximation, in which the
leading-order geometry is described by a Tolman--Bondi dust
solution~\cite{Tolman1934,Bondi1947}.

By Corollary~\ref{cor:small-initial-compactness}, the dimensionless product $\mu r_0$ in our construction can take any
prescribed positive value; recall that $r_0$ is the area radius of the
central FLRW ball.
Our construction imposes no smallness or largeness assumption on $\mu r_0$.

\subsection{Gluing methods and related constructions}\label{sec:proof-strategy}

A main technique in this paper is to glue an interior region to an
exterior region while satisfying the constraint equations. This approach
goes back to Corvino \cite{Corvino2000} and was further developed by
Corvino--Schoen \cite{CorvinoSchoen2006}. Their work established methods
for gluing suitable vacuum initial data to Schwarzschild or Kerr
initial data in an exterior region. In a series of papers
\cite{AretakisCzimekRodnianski2023,AretakisCzimekRodnianski2024,AretakisCzimekRodnianski2025},
Aretakis--Czimek--Rodnianski developed a characteristic gluing method
for the Einstein vacuum equations and used it to glue strongly
asymptotically flat characteristic data to data of a suitably chosen Kerr spacetime.

Gluing methods have been applied to construct asymptotically flat
initial data for trapped surface formation. In vacuum, Li--Yu
\cite{LiYu2015} combined Christodoulou's short-pulse method
with the Corvino--Schoen gluing method to obtain such data. Using characteristic gluing, Kehle--Unger \cite{KehleUnger2024} constructed asymptotically flat
vacuum Cauchy data that collapse to very slowly
rotating Kerr black holes. In \cite{KehleUngerCharged}, they also constructed spherically symmetric initial data that collapse to Reissner--Nordstr\"om black holes for the Einstein--Maxwell--charged scalar field system.
Giorgi--Shen--Wan \cite{GiorgiShenWan2025} combined short-pulse constructions with spacelike gluing to obtain
initial data whose evolution forms several causally independent trapped surfaces for the Einstein--Maxwell--Klein--Gordon system.

Our construction uses a contracting FLRW ball as the interior region. The aim is to extend the initial data on the interior ball to complete initial data that are untrapped and
asymptotically flat. We glue the interior data to Schwarzschild exterior
data through a transition region. The data in the transition region satisfy the constraint
equations. In the spherically
symmetric Einstein--Klein--Gordon system, the constraint equations can
be reduced to an ordinary differential equation for the Hawking mass.
We analyze this ODE explicitly and choose the transition data so that
every centered sphere has positive outward null expansion. This
excludes weakly outer trapped surfaces on the initial
hypersurface. Since the initial data in the interior ball are unchanged, the resulting
spacetime agrees with the FLRW solution in the future domain of
dependence of the interior ball, where trapped surfaces form.

In fact, the idea of gluing a homogeneous collapsing interior to a vacuum exterior has a long history. Oppenheimer--Snyder \cite{OppenheimerSnyder1939} constructed a model of black hole formation by matching a homogeneous dust interior to a Schwarzschild exterior. In this model, the matter density is discontinuous at the boundary of the dust cloud. More recently, Andr\'easson--Rein \cite{AndreassonRein2025} established an Oppenheimer--Snyder type collapse result for the massive Einstein--Vlasov system. Their initial data contain a spatially homogeneous interior and a regular transition region to a vacuum exterior. For a suitable family approximating the Oppenheimer--Snyder dust model, they proved that trapped surfaces form from untrapped initial data. They approximated Oppenheimer--Snyder dust data by regular Einstein--Vlasov data with a small velocity dispersion, using a transition layer to connect a spatially homogeneous central core to an exterior vacuum region. They showed that the evolution in a shrinking central region stays close to the Oppenheimer--Snyder dust solution long enough for a trapped surface to form. Note that the homogeneous dust interior in the original Oppenheimer--Snyder construction is spatially flat on each slice of constant comoving time. Cohen \cite{Cohen1968} extended this construction by considering interior regions with nonzero spatial curvature.

\subsection{Outline}

\Cref{sec:constraints} constructs the asymptotically flat extension for the interior FLRW ball,
establishes the compactness estimate, and excludes weakly outer
trapped surfaces for the completed initial data. \Cref{sec:geometry} studies the FLRW solution and
derives sufficient conditions for trapped surface formation and curvature
blow-up inside the future domain of dependence of the interior FLRW ball. \Cref{sec:comparison} compares
the radius thresholds and proves the results stated in
\cref{sec:main-results}. \Cref{sec:discussion} discusses the stability of trapped surface
formation under perturbations of the initial data and questions
for further study.

\subsection{Acknowledgements}

The author is greatly indebted to his advisor Shiwu Yang for suggesting and leading to the problem and for helpful discussions.

\section{Untrapped asymptotically flat initial data}\label{sec:constraints}

In this section, we complete the initial data on the interior FLRW ball
by gluing them to an exterior Schwarzschild slice. For a spacelike
hypersurface $\Sigma$ with induced metric $h$ and future unit normal $n$,
set
\begin{equation}\label{eq:data-convention}
 K=\frac12\mathcal L_nh,
 \qquad \Phi=\phi|_\Sigma,
 \qquad \Pi=n(\phi)|_\Sigma.
\end{equation}
The Hamiltonian and momentum constraint equations are
\begin{align}
 \Scal_h+(\tr_hK)^2-|K|_h^2
 &=2\bigl(\Pi^2+|D\Phi|_h^2+\mu^2\Phi^2\bigr),
 \label{eq:hamiltonian}\\*
 D^jK_{ij}-D_i\tr_hK&=2\Pi D_i\Phi,
 \label{eq:momentum}
\end{align}
where $D$ is the Levi--Civita connection of $h$. Smooth data satisfying
these constraint equations have a maximal globally hyperbolic
development, unique up to isometry
\cite{ChoquetBruhatGeroch1969,ChoquetBruhat2009}.

In spherical symmetry, write
\begin{equation}\label{eq:areal-data}
 h=B(r)^2\dd r^2+r^2\dd\omega^2,
 \qquad
 K^i{}_j=\operatorname{diag}\bigl(q(r),v(r),v(r)\bigr),
\end{equation}
where $B>0$ and the entries are ordered according to $(r,\theta,\varphi)$.
The scalar data are radial functions, $\Phi=\Phi(r)$ and $\Pi=\Pi(r)$.
Define the Riemannian Hawking mass $\msp$, the spacetime Hawking
mass $\HM$, and the function $\Delta$ by
\begin{equation}\label{eq:masses}
 B^{-2}=1-\frac{2\msp}{r},
 \qquad
 \HM=\msp+\frac12r^3v^2,
 \qquad
 \Delta=1-\frac{2\HM}{r}=B^{-2}-r^2v^2.
\end{equation}
The condition $\Delta>0$ ensures that the corresponding sphere at
constant $r$ is untrapped. We can now express the constraint equations
in terms of the newly introduced quantities $q$, $v$, and $\HM$.

\begin{lemma}\label{lem:reduced-constraints}
For smooth data of the form \eqref{eq:areal-data}, the constraint
equations \eqref{eq:hamiltonian}--\eqref{eq:momentum} are equivalent to
\begin{align}
 q&=v+rv'+r\Pi\Phi',
 \label{eq:q-constraint}\\
 \HM'&=\frac{r^2}{2}\left[
 (\Pi-rv\Phi')^2+\Delta(\Phi')^2+\mu^2\Phi^2
 \right],
 \label{eq:M-constraint}
\end{align}
where a prime denotes differentiation with respect to $r$.
Equivalently, the mass satisfies the linear ordinary differential
equation
\begin{equation}\label{eq:M-linear}
 \HM'+r(\Phi')^2\HM
 =\frac{r^2}{2}\left[
 (\Pi-rv\Phi')^2+(\Phi')^2+\mu^2\Phi^2
 \right].
\end{equation}
In particular, \eqref{eq:M-constraint} gives $\HM'\ge0$ wherever
$\Delta>0$.
\end{lemma}

\begin{proof}
Direct calculation in the metric \eqref{eq:areal-data} gives
\begin{align*}
 \Scal_h&=\frac{4\msp'}{r^2},
 &(\tr_hK)^2&=(q+2v)^2,
 &|K|_h^2&=q^2+2v^2,\\
 |D\Phi|_h^2&=\frac{(\Phi')^2}{B^2},
 &D^jK_{rj}&=q'+\frac2r(q-v),
 &D_r\tr_hK&=q'+2v',
\end{align*}
and
\begin{gather*}
 D^jK_{\theta j}=D^jK_{\varphi j}
 =D_\theta\tr_hK=D_\varphi\tr_hK=0,\\
 \Pi D_r\Phi=\Pi\Phi',
 \qquad \Pi D_\theta\Phi=\Pi D_\varphi\Phi=0.
\end{gather*}
Thus the radial component of the momentum constraint
\eqref{eq:momentum} gives \eqref{eq:q-constraint}, while its angular
components vanish identically. The Hamiltonian constraint becomes
\[
 \frac{2\msp'}{r^2}+2qv+v^2
 =\Pi^2+\frac{(\Phi')^2}{B^2}+\mu^2\Phi^2.
\]
Substituting \eqref{eq:q-constraint} for $q$ and using
$\msp=\HM-\frac12r^3v^2$ from \eqref{eq:masses} gives
\[
 \HM'=\frac{r^2}{2}
 \left[\Pi^2-2rv\Pi\Phi'+\frac{(\Phi')^2}{B^2}+\mu^2\Phi^2\right].
\]
Substituting $B^{-2}=\Delta+r^2v^2$ and completing the square yields
\[
 \HM'=\frac{r^2}{2}
 \left[(\Pi-rv\Phi')^2+\Delta(\Phi')^2+\mu^2\Phi^2\right],
\]
which is \eqref{eq:M-constraint}. Replacing $\Delta$ by $1-2\HM/r$
gives \eqref{eq:M-linear}. These steps are reversible.
\end{proof}

By \eqref{eq:completed-data-intro}, on the interior FLRW ball we have
\begin{equation}\label{eq:ball-areal-variables}
 B^{-2}=1-\kappa r^2,\qquad q=v=H_0,\qquad
 \Phi=\Phi_0,\qquad \Pi=\Pi_0
 \quad(0\le r<r_0).
\end{equation}
Using \eqref{eq:masses} and \eqref{eq:E0-H0}, we obtain
\begin{equation}\label{eq:ball-masses}
 \msp=\frac12\kappa r^3,\qquad
 \HM=\frac12(\kappa+H_0^2)r^3=\frac{E_0r^3}{6},\qquad
 \Delta=1-\frac{E_0r^2}{3},
\end{equation}
where
\begin{equation}\label{eq:initial-scalar-energy}
 E_0:=\Pi_0^2+\mu^2\Phi_0^2=3(\kappa+H_0^2).
\end{equation}
$H_0$ is defined by \eqref{eq:E0-H0}. Note that the Hawking mass $\HM$ and $\Delta$ on the interior ball
are independent of $\kappa$ when the scalar data $\Phi_0$ and $\Pi_0$
are fixed.

The following proposition gives a sufficient condition for the
existence of an untrapped, asymptotically flat completion of the
interior FLRW ball.

\begin{proposition}\label{prop:completion}
Fix $\mu>0$, $(\Phi_0,\Pi_0)\ne(0,0)$, and $\kappa<E_0/3$, with
$E_0$ defined by \eqref{eq:initial-scalar-energy}. If
\begin{equation}\label{eq:completion-condition}
 0<r_0<\Rglue
 :=\min\left\{\sqrt{\frac3{E_0}},\frac1{\mu|\Phi_0|}\right\},
\end{equation}
then the data \eqref{eq:ball-areal-variables} on the interior FLRW
ball extend smoothly to complete, spherically symmetric,
asymptotically flat data on $\R^3$. The extension satisfies the
constraint equations, has $\Delta>0$ everywhere, and coincides with
Schwarzschild data \eqref{eq:completed-data-intro} for $r>R_0$, for some constants $M>0$ and $R_0>r_0$.
The second entry in the minimum is understood as $+\infty$ when
$\Phi_0=0$.
\end{proposition}

\begin{proof}
Choose $r_a$, $r_b$, and $R_0$ such that
$r_0<r_a<r_b<R_0<\Rglue$. We prescribe $\Phi$, $\Pi$, and $v$,
solve \eqref{eq:M-linear} for $\HM$ with
$\HM(r_0)=E_0r_0^3/6$, and recover $q$ and $B$ from
\eqref{eq:q-constraint} and $B^{-2}=\Delta+r^2v^2$, where
$\Delta=1-2\HM/r$. All cutoffs below are constant near their endpoints.

On $[r_0,r_a]$, set $\Phi=\Phi_0$ and $v=H_0$, and cut off $\Pi$
smoothly from $\Pi_0$ to zero with $|\Pi|\le|\Pi_0|$.
Since $\Phi'=0$, \eqref{eq:M-linear} gives
\begin{equation}\label{eq:zone-one-mass}
 \HM'=\frac{r^2}{2}(\Pi^2+\mu^2\Phi_0^2)
 \le\frac{E_0r^2}{2}.
\end{equation}
Integrating \eqref{eq:zone-one-mass} and using
$\HM(r_0)=E_0r_0^3/6$, we obtain
\begin{equation}\label{eq:zone-one-mass-bound}
 \HM(r)\le\frac{E_0r^3}{6}.
\end{equation}
It follows that
\begin{equation}\label{eq:zone-one-delta}
 \Delta(r)\ge1-\frac{E_0r^2}{3}
 >1-\frac{E_0\Rglue^2}{3}\ge0.
\end{equation}

On $[r_a,r_b]$, cut off $\Phi$ smoothly from $\Phi_0$ to zero with
$|\Phi|\le|\Phi_0|$, and set
\begin{equation}\label{eq:zone-two-profiles}
 v=H_0,\qquad \Pi=rv\Phi'.
\end{equation}
Substituting \eqref{eq:zone-two-profiles} and
$\HM=\frac r2(1-\Delta)$ into \eqref{eq:M-constraint} gives
\begin{equation}\label{eq:delta-equation}
 r\Delta'+\bigl(1+r^2(\Phi')^2\bigr)\Delta
 =1-\mu^2r^2\Phi^2.
\end{equation}
Since $r_b<\Rglue$ and $|\Phi|\le|\Phi_0|$, its right-hand side
satisfies
\[
 1-\mu^2r^2\Phi^2
 \ge1-\mu^2r_b^2\Phi_0^2>0.
\]
If $\Delta$ vanished at some $r\in(r_a,r_b]$,
\eqref{eq:delta-equation} would give $\Delta'(r)>0$.
Such a point cannot be the first zero reached from positive values.
Since $\Delta(r_a)>0$ by \eqref{eq:zone-one-delta}, we conclude that
$\Delta>0$ throughout $[r_a,r_b]$.

On $[r_b,R_0]$, keep $\Phi=\Pi=0$ and cut off $v$ smoothly from
$H_0$ to zero. Equation \eqref{eq:M-linear} gives $\HM'=0$, so
$\HM=\HM(r_b)$ is constant. By \eqref{eq:M-constraint} and the
positivity of $\Delta$ on the preceding annuli,
$\HM(r_b)\ge\HM(r_0)>0$. Thus $\Delta(r)=1-2\HM(r_b)/r$
increases on $[r_b,R_0]$ and remains positive.

For $r\ge R_0$, keep $\Phi=\Pi=v=0$. Then $q=0$ and $\HM=M$,
where $M:=\HM(R_0)>0$. The Cauchy data in this exterior region
coincide with the Schwarzschild data
\begin{equation}\label{eq:schwarzschild-end}
 h=\frac{\dd r^2}{1-2M/r}+r^2\dd\omega^2,
 \qquad \Phi=\Pi=0,\qquad K=0.
\end{equation}
The construction keeps $\Phi$, $\Pi$, $v$, $q$, $\HM$, $\Delta$, and
$B$ smooth, with $\Delta>0$, throughout $r>0$. The resulting initial
data coincide with the FLRW data for $r<r_0$ and with the
Schwarzschild data for $r>R_0$. The regular center and the
Schwarzschild exterior give a complete, asymptotically flat
hypersurface diffeomorphic to $\R^3$.
\end{proof}

For the small-compactness conclusion in
Corollary~\ref{cor:small-initial-compactness}, it is convenient to choose
$\Phi_0=0$. Then only the cutoffs of $\Pi$ and $v$ need to be
considered. The following proposition treats this case explicitly.

\begin{proposition}\label{prop:small-initial-mass-ratio}
Fix $\mu>0$ and $\Pi_0\ne0$, set $\Phi_0=0$, and let
$\kappa<\Pi_0^2/3$. Choose $r_0$ and $R_0$ such that
\[
 0<r_0<R_0<\frac{\sqrt3}{|\Pi_0|}=\Rglue.
\]
Then the interior FLRW data \eqref{eq:ball-areal-variables} extend
smoothly to untrapped, asymptotically flat initial data on a complete
Cauchy hypersurface $\Sigma\cong\R^3$. The data coincide with
Schwarzschild data for $r>R_0$ and satisfy
\begin{equation}\label{eq:general-initial-mass-ratio}
 0<\sup_\Sigma\frac{2\HM}{r}
 \le\frac{\Pi_0^2R_0^2}{3}=\left(\frac{R_0}{\Rglue}\right)^2<1.
\end{equation}
\end{proposition}

\begin{proof}
Extend the interior FLRW ball as in the proof of
Proposition~\ref{prop:completion}, choosing $r_a$ such that $r_0<r_a<R_0$.
Since $\Phi_0=0$, we can set $\Phi=0$ throughout $\Sigma$ and omit
the cutoff of $\Phi$. On $[r_0,r_a]$, keep $v=H_0$ and cut off
$\Pi$ smoothly from $\Pi_0$ to zero with $|\Pi|\le|\Pi_0|$.
Then keep $\Pi=0$ for $r\ge r_a$, cut off $v$ smoothly to zero on
$[r_a,R_0]$, and set $v=0$ for $r\ge R_0$.
By \eqref{eq:zone-one-mass-bound} and the interior mass formula
\eqref{eq:ball-masses},
\[
 0<\frac{2\HM(r)}r
 \le\frac{\Pi_0^2r^2}{3}
 \le\frac{\Pi_0^2r_a^2}{3}
 <\frac{\Pi_0^2R_0^2}{3}<1
 \qquad(0<r\le r_a).
\]
For $r\ge r_a$, the mass $\HM(r)=\HM(r_a)$ is constant, so
\[
 0<\frac{2\HM(r)}r
 =\frac{2\HM(r_a)}r
 \le\frac{2\HM(r_a)}{r_a}
 \le\frac{\Pi_0^2r_a^2}{3}
 <\frac{\Pi_0^2R_0^2}{3}<1.
\]
These inequalities give \eqref{eq:general-initial-mass-ratio}.
The data are smooth and complete by the same argument as in
Proposition~\ref{prop:completion}. Since $\Delta>0$,
Lemma~\ref{lem:touching} below shows that they are untrapped.
\end{proof}

In \cref{sec:comparison}, we show that, by taking $\kappa$ sufficiently
negative, we can choose $r_0$ and $R_0$ with
$\Rtrap<r_0<R_0<\Rglue$ while making $R_0/\Rglue$ arbitrarily small.
It follows from
\[
 \sup_\Sigma\frac{2\HM}{r}
 \le\frac{\Pi_0^2R_0^2}{3}
 =\left(\frac{R_0}{\Rglue}\right)^2
\]
that the supremum can be made arbitrarily small. This explains the small-compactness conclusion
of Corollary~\ref{cor:small-initial-compactness}; the case of prescribed $r_0$
is treated in its proof.

The quantity $\Delta=1-2\HM/r$ determines whether the centered
spheres are untrapped. For $r>0$, let $S_r\subset\Sigma$ denote the
centered sphere of area radius $r$, whose area is $4\pi r^2$.
Its future outgoing and ingoing null normals are
$\ell_\pm=n\pm B^{-1}\partial_r$. For a local orthonormal tangent
frame $e_1,e_2$ on $S_r$, the corresponding null expansions are defined by
\begin{equation}\label{eq:null-expansion-definition}
 \theta_\pm(S_r)=\sum_{a=1}^2g(\nabla^g_{e_a}\ell_\pm,e_a),
\end{equation}
where $\nabla^g$ is the Levi--Civita connection of $g$.
For the data \eqref{eq:areal-data}, this gives
\begin{equation}\label{eq:initial-centered-expansions}
 \theta_\pm(S_r)=\frac2r\bigl(rv\pm B^{-1}\bigr).
\end{equation}
By \eqref{eq:masses} and \eqref{eq:initial-centered-expansions},
\begin{equation}\label{eq:delta-expansions}
 \Delta=B^{-2}-r^2v^2
 =-\frac{r^2}{4}\theta_+(S_r)\theta_-(S_r).
\end{equation}
Since $\theta_+>\theta_-$, the condition $\Delta>0$ implies
$\theta_+>0>\theta_-$. The following lemma shows that if $\Delta>0$
throughout the initial hypersurface, then it contains no closed
weakly outer trapped surface. Here surfaces are smooth, connected,
embedded, and two-sided, with the outward normal pointing toward
the asymptotically flat end.

\begin{lemma}\label{lem:touching}
Let $h$ and $K$ be spherically symmetric and have the form
\eqref{eq:areal-data} on the initial hypersurface $\Sigma\cong\R^3$.
If every centered sphere $S_r$ has positive outward future null
expansion $\theta_+(S_r)>0$, then $\Sigma$ contains no smooth closed
weakly outer trapped surface.
\end{lemma}

\begin{proof}
Suppose that a smooth closed weakly outer trapped surface $\Gamma$
exists. Let $R=\max_\Gamma r$, attained at $p\in\Gamma$.
The sphere $S_R=\{r=R\}$ encloses $\Gamma$ and touches it at $p$.
At this point they have the same tangent plane and outward unit
normal $\nu=B(R)^{-1}\partial_r$. In particular, $\nu(r)>0$.

Since $p$ is a maximum point of $r|_\Gamma$, the Laplacian of $r$
on $\Gamma$ satisfies $\Delta_\Gamma r(p)\le0$. Since $r$ is
constant on $S_R$, we have $\Delta_{S_R}r=0$.
The formula for the Laplace--Beltrami operator on a surface gives,
at $p$,
\begin{align*}
 \Delta_\Gamma r
 &=\tr_{T_p\Gamma}(\operatorname{Hess}_h r)-H_\Gamma\nu(r)\le0,\\
 \Delta_{S_R}r
 &=\tr_{T_pS_R}(\operatorname{Hess}_h r)-H_{S_R}\nu(r)=0.
\end{align*}
Here $\operatorname{Hess}_h r$ is the Hessian of $r$ on $(\Sigma,h)$,
and $H_\Gamma$ and $H_{S_R}$ are the mean curvatures with respect to
the outward normals, with the convention $H=\operatorname{div}\nolimits_S\nu$.
Because $T_p\Gamma=T_pS_R$, the two traces of
$\operatorname{Hess}_h r$ agree. The preceding identities and
$\nu(r)>0$ therefore imply $H_\Gamma(p)\ge H_{S_R}(p)$.

The outward future null normal to $\Gamma$ is $\ell_+=n+\nu$.
For a local orthonormal tangent frame $e_1,e_2$ on $\Gamma$, the
definition \eqref{eq:null-expansion-definition} of the null expansion gives
\begin{align*}
 \theta_+(\Gamma)
 &=\sum_{a=1}^2g(\nabla^g_{e_a}\ell_+,e_a)\\
 &=\sum_{a=1}^2g(\nabla^g_{e_a}n,e_a)
   +\sum_{a=1}^2g(\nabla^g_{e_a}\nu,e_a)\\
 &=\sum_{a=1}^2K(e_a,e_a)
   +\sum_{a=1}^2h(D_{e_a}\nu,e_a)
 =\tr_\Gamma K+H_\Gamma.
\end{align*}
We used the convention \eqref{eq:data-convention} and the fact that
$D_{e_a}\nu$ is the component of $\nabla^g_{e_a}\nu$ tangent to
$\Sigma$. At $p$, the common tangent plane also gives
$\tr_{T_p\Gamma}K=\tr_{T_pS_R}K$. Hence
\[
 \theta_+(\Gamma)(p)
 =H_\Gamma(p)+\tr_{T_p\Gamma}K
 \ge H_{S_R}(p)+\tr_{T_pS_R}K
 =\theta_+(S_R)(p)>0.
\]
This contradicts the assumption that $\Gamma$ is weakly outer trapped.
\end{proof}

Lemma~\ref{lem:touching} also follows from the maximum principle in
\cite[Proposition~3.1]{AshtekarGalloway2005}. Together with
\eqref{eq:initial-centered-expansions}, it shows that the data
constructed in Proposition~\ref{prop:completion} contain no closed weakly outer
trapped surface and are therefore initially untrapped.

\section{Formation of trapped surfaces and curvature singularities}\label{sec:geometry}

We study the evolution inside the future domain of dependence of the
initial FLRW ball. For certain scalar data $(\Phi_0,\Pi_0)$, there is a
range of initial area radii for which the ball admits an untrapped
completion as in Proposition~\ref{prop:completion} and a trapped surface forms in
its future domain of dependence.

The metric in the FLRW region is given by \eqref{eq:flrw-intro}, and its
evolution is governed by \eqref{eq:flrw-constraint-intro}--\eqref{eq:flrw-kg-intro}.
Recall that the centered sphere at $(t,s)$ has area radius
\begin{equation}\label{eq:area-radius}
 r(t,s)=a(t)S_\kappa(s),
\end{equation}
where, as in \eqref{eq:S-kappa-intro},
\[
 S_\kappa(s)=
 \begin{cases}
  \sin(\sqrt{\kappa}\,s)/\sqrt{\kappa}, & \kappa>0,\\
  s, & \kappa=0,\\
  \sinh(\sqrt{-\kappa}\,s)/\sqrt{-\kappa}, & \kappa<0.
 \end{cases}
\]
The function $S_\kappa$ satisfies
\begin{equation}\label{eq:curvature-identities}
 (S_\kappa')^2+\kappa S_\kappa^2=1.
\end{equation}
We work on its increasing branch. The outward and inward
future radial null vectors are
\begin{equation}\label{eq:null-vectors}
 L_\pm=\partial_t\pm a^{-1}\partial_s.
\end{equation}
By \eqref{eq:area-radius} and \eqref{eq:null-vectors}, the null
expansions of the sphere of area $4\pi r(t,s)^2$ are
\begin{equation}\label{eq:flrw-expansions}
 \begin{aligned}
 \theta_\pm(t,s)
 &=\frac{L_\pm(4\pi r(t,s)^2)}{4\pi r(t,s)^2}
 =\frac{2L_\pm r(t,s)}{r(t,s)}=2H(t)\pm\frac{2S_\kappa'(s)}{a(t)S_\kappa(s)},
 \end{aligned}
\end{equation}
where $H(t)=\dot a(t)/a(t)$.

Our choice $H(0)=-\sqrt{(\Pi_0^2+\mu^2\Phi_0^2)/3-\kappa}<0$ ensures
that the solution is initially contracting. Recall that $[0,T^*)$ is
the maximal interval starting at $t=0$ on which the solution is regular,
$a(t)>0$, and $H(t)<0$. On this interval, $H(t)<0$ and $S_\kappa'(s)>0$,
so the inward expansion $\theta_-$ is negative for $s>0$.
The outward expansion is negative precisely when
$|\dot a|S_\kappa(s)>S_\kappa'(s)$. Squaring both sides and using
\eqref{eq:curvature-identities}, we obtain
\begin{equation}\label{eq:trapped-equivalences}
 \begin{aligned}
 \theta_+(t,s)<0
 &\quad\Longleftrightarrow\quad S_\kappa(s)^2\bigl(\dot a(t)^2+\kappa\bigr)>1\\
 &\quad\Longleftrightarrow\quad S_\kappa(s)>\frac{1}{\sqrt{\dot a(t)^2+\kappa}}\\
 &\quad\Longleftrightarrow\quad r(t,s)>\frac{a(t)}{\sqrt{\dot a(t)^2+\kappa}}.
 \end{aligned}
\end{equation}
The quantity on the right is the area radius of the apparent horizon on
the slice at time $t$ of the FLRW solution; see also \cite{MisnerSharp1964,Faraoni2011}.
By \eqref{eq:flrw-constraint-intro},
\begin{equation}\label{eq:apparent-area-radius}
 r_{\mathrm{AH}}(t):=\frac{a(t)}{\sqrt{\dot a(t)^2+\kappa}}
 =\frac{\sqrt{3}}{\sqrt{\dot\phi(t)^2+\mu^2\phi(t)^2}}.
\end{equation}
If $\phi$ and $\dot\phi$ both vanished at a regular time, uniqueness
for the Klein--Gordon equation would give $\phi\equiv0$, contrary to
$(\Phi_0,\Pi_0)\ne(0,0)$. Thus the denominators in
\eqref{eq:apparent-area-radius} are positive.

The corresponding geodesic coordinate of the apparent horizon is
\begin{equation}\label{eq:horizon-coordinate}
 s_{\mathrm{AH}}(t):=S_\kappa^{-1}\left(\frac{r_{\mathrm{AH}}(t)}{a(t)}\right)
 =S_\kappa^{-1}\left(\frac{1}{\sqrt{\dot a(t)^2+\kappa}}\right),
\end{equation}
where the inverse on the increasing branch is given by
\begin{equation}\label{eq:inverse-radius}
 S_\kappa^{-1}(x)=
 \begin{cases}
  \kappa^{-1/2}\arcsin(\sqrt{\kappa}\,x), & \kappa>0,\\
  x, & \kappa=0,\\
  (-\kappa)^{-1/2}\operatorname{arsinh}(\sqrt{-\kappa}\,x), & \kappa<0.
 \end{cases}
\end{equation}
Thus $s_{\mathrm{AH}}(t)$ is well defined throughout the contracting
interval. In particular, if $\kappa>0$, then $\dot a(t)<0$ gives
$(\dot a(t)^2+\kappa)^{-1/2}<\kappa^{-1/2}$, so
\[
 0<s_{\mathrm{AH}}(t)<\frac{\pi}{2\sqrt{\kappa}}.
\]

We can now define the trapping threshold. Set
\begin{equation}\label{eq:trapping-coordinates}
 \begin{aligned}
 \eta(t)&:=\int_0^t\frac{\dd\tau}{a(\tau)},\\
 \widetilde{s}(t)&:=\eta(t)+s_{\mathrm{AH}}(t),\qquad
 s_{\mathrm{trap}}:=\inf_{0<t<T^*}\widetilde{s}(t).
 \end{aligned}
\end{equation}
The threshold in initial area radius is
\begin{equation}\label{eq:trapping-radius}
 \Rtrap=S_\kappa(s_{\mathrm{trap}}),
\end{equation}
in agreement with \eqref{eq:thresholds-intro}. Here $\eta(t)$ is the
comoving distance traveled by a radial light ray from time $0$ to time $t$.

Let the initial homogeneous ball and its area radius be
\begin{equation}\label{eq:U-geodesic}
 \U=\{0\le s<s_0\},\qquad r_0=S_\kappa(s_0).
\end{equation}
For $\kappa>0$, we take $s_0<\pi/(2\sqrt{\kappa})$.
By finite propagation speed and geometric uniqueness for the
Einstein--Klein--Gordon system, the solution in the future domain of
dependence $D^+(\U)$ agrees with the corresponding FLRW solution up to
isometry \cite{ChoquetBruhatGeroch1969,ChoquetBruhat2009}.
Here we restrict this identification to the contracting interval
$0\le t<T^*$. The following proposition gives a criterion for trapped surface
formation in this region.

\begin{proposition}\label{prop:causal-criterion}
The future domain of dependence $D^+(\U)$ contains a trapped centered
sphere at some time before $T^*$ if $s_0>s_{\mathrm{trap}}$, or equivalently,
if
\begin{equation}\label{eq:causal-criterion}
 r_0>\Rtrap.
\end{equation}
\end{proposition}

\begin{proof}
We use FLRW coordinates on the region identified with $D^+(\U)$.
The inward radial null ray starting at $(t,s)=(0,s_0)$ has coordinate
\begin{equation}\label{eq:inward-ray}
 s(t)=s_0-\int_0^t\frac{\dd\tau}{a(\tau)}=s_0-\eta(t)
\end{equation}
until it reaches the center. Since every causal curve satisfies
$|\dd s/\dd t|\le a(t)^{-1}$, for $0<t<T^*$ with $\eta(t)<s_0$ the
part of the FLRW slice $\Sigma_t$ in the interior of $D^+(\U)$ is
\begin{equation}\label{eq:Dplus-slice}
 \{p\in\Sigma_t:0\le s(p)<s_0-\eta(t)\}.
\end{equation}
If
\begin{equation}\label{eq:trapped-sphere-room}
 s_{\mathrm{AH}}(t)<s_0-\eta(t),
\end{equation}
we can choose $s$ such that $s_{\mathrm{AH}}(t)<s<s_0-\eta(t)$.
When $\kappa>0$, this choice also gives $s<\pi/(2\sqrt{\kappa})$.
The centered sphere at $(t,s)$ then lies in $D^+(\U)$ and is trapped
by \eqref{eq:trapped-equivalences} and \eqref{eq:horizon-coordinate}.

If $s_0>s_{\mathrm{trap}}$, the definition of the infimum gives a time
$0<t<T^*$ such that $s_0>\eta(t)+s_{\mathrm{AH}}(t)$. Thus
\eqref{eq:trapped-sphere-room} holds. Finally, $S_\kappa$ is strictly
increasing on the branch under consideration. For $\kappa>0$, the
threshold also lies on this branch, since
$s_{\mathrm{trap}}\le\widetilde{s}(0)=s_{\mathrm{AH}}(0)<\pi/(2\sqrt{\kappa})$
by continuity at $t=0$. Hence $s_0>s_{\mathrm{trap}}$ is equivalent
to \eqref{eq:causal-criterion}.
\end{proof}

Besides trapped surface formation, we also consider curvature singularity
formation along the central geodesic in $D^+(\U)$. The following
proposition gives a sufficient condition.

\begin{proposition}\label{prop:causal-collapse}
Let $[0,T^*)$ be the maximal regular contracting interval of the FLRW
solution. Suppose that
\[
 T^*<\infty,\qquad a^{-1}\in L^1(0,T^*),\qquad
 a(t)\to0\quad(t\to T^*).
\]
Let the initial Cauchy data agree with the FLRW data on the interior
ball $\U=\{0\le s<s_0\}$. If
\begin{equation}\label{eq:causal-collapse}
 s_0>\eta_*:=\int_0^{T^*}\frac{\dd t}{a(t)},
\end{equation}
then the central comoving geodesic $\gamma$, given by $s=0$, remains in
$\operatorname{int}D^+(\U)$ for $0<t<T^*$ and satisfies
\begin{equation}\label{eq:central-curvature-blowup}
 \operatorname{Scal}_g(\gamma(t))\longrightarrow-\infty
 \qquad(t\to T^*).
\end{equation}
\end{proposition}

\begin{proof}
Set
\begin{equation}\label{eq:scalar-energy}
 E(t):=\dot\phi(t)^2+\mu^2\phi(t)^2.
\end{equation}
The Klein--Gordon equation
\eqref{eq:flrw-kg-intro} gives
\begin{equation}\label{eq:contracting-energy}
 \dot E(t)=-6H(t)\dot\phi(t)^2\ge0,
\end{equation}
since $H(t)<0$ on $[0,T^*)$. If $E$ were bounded on this interval,
the Friedmann equation \eqref{eq:flrw-constraint-intro} would give
\[
 \dot a(t)^2=\frac{a(t)^2E(t)}{3}-\kappa.
\]
Since $0<a(t)\le a(0)=1$, this would imply $|\dot a(t)|\le C$ for
some constant $C$. Consequently,
\[
 \log\frac{a(0)}{a(t)}
 =\int_0^t\frac{-\dot a(\tau)}{a(\tau)}\,\dd\tau
 \le C\eta_*,
\]
contrary to $a(t)\to0$. Thus $E$ is unbounded and, being nondecreasing,
satisfies $E(t)\to\infty$ as $t\to T^*$.

The definition of $E$ gives $|\dot\phi|\le\sqrt E$. For any
$0\le t_0<t<T^*$, the monotonicity of $E$ therefore yields
\[
 \begin{aligned}
 |\phi(t)|
 &\le |\phi(t_0)|+\int_{t_0}^t\sqrt{E(\tau)}\,\dd\tau\\
 &\le |\phi(t_0)|+(t-t_0)\sqrt{E(t)}.
 \end{aligned}
\]
Dividing by $\sqrt{E(t)}$ and letting $t\to T^*$, we obtain
\[
 0\le\limsup_{t\to T^*}\frac{|\phi(t)|}{\sqrt{E(t)}}
 \le T^*-t_0.
\]
Letting $t_0\to T^*$ gives $\phi(t)^2/E(t)\to0$.
Taking the trace of the Einstein equation \eqref{eq:ekg-einstein}
and using $\phi=\phi(t)$ gives the scalar curvature of the FLRW solution:
\begin{equation}\label{eq:flrw-scalar-curvature}
 \begin{aligned}
 \operatorname{Scal}_{g_{\mathrm{FLRW}}}
 &=2g^{\alpha\beta}
   \partial_\alpha\phi\,\partial_\beta\phi+4\mu^2\phi^2\\
 &=4\mu^2\phi^2-2\dot\phi^2
 =6\mu^2\phi^2-2E\\
 &=E\left(-2+6\mu^2\frac{\phi^2}{E}\right)
 \longrightarrow-\infty
 \end{aligned}
\end{equation}
as $t\to T^*$. Finally, $\eta(t)<\eta_*<s_0$ for every $t<T^*$. By
\eqref{eq:Dplus-slice}, the central comoving geodesic remains in the
interior of $D^+(\U)$, where the spacetime agrees with the FLRW
solution. The scalar curvature therefore also blows up along $\gamma$, proving
\eqref{eq:central-curvature-blowup}.
\end{proof}

The condition $a(t)\to0$ at a finite time does not by itself imply
curvature blow-up. For comparison, the vacuum solution
\[
 \kappa=-1,\qquad \phi=0,\qquad a(t)=1-t
 \quad(0\le t<1)
\]
is a contracting Milne region of Minkowski spacetime. Although $a(t)\to0$ as $t\to1$, the spacetime curvature
vanishes.
Note that $\int_0^1a(t)^{-1}\dd t=\infty$.  Such examples also exist for nonzero scalar data when
$\kappa<0$: the homogeneous massive Einstein--Klein--Gordon equations
\eqref{eq:flrw-constraint-intro}--\eqref{eq:flrw-kg-intro}
admit special contracting solutions with
\[
 a(t)\sim\sqrt{-\kappa}(T^*-t),
 \qquad \phi(t)\to\phi_*\ne0,
 \qquad \dot\phi(t)\to0
 \quad(t\to T^*).
\]
Their spacetime curvature remains
bounded, while $a^{-1}$ is not integrable on $[0, T^*)$. These examples illustrate that $a(t)\to0$ alone does not imply curvature
blow-up.
The integrability of $a^{-1}$ on $[0, T^*)$ guarantees curvature blow-up under the other hypotheses of Proposition~\ref{prop:causal-collapse}.

For $\kappa\le0$, $S_\kappa$ is strictly increasing. Thus, when
$\eta_*<\infty$, condition \eqref{eq:causal-collapse} is equivalent to
\[
 r_0=S_\kappa(s_0)>S_\kappa(\eta_*)=\Rcoll.
\]
This gives the following corollary.

\begin{corollary}\label{cor:causal-collapse}
Let $\kappa\le0$, and let $[0,T^*)$ be the maximal regular contracting interval of the FLRW
solution. Suppose that
\[
 T^*<\infty,\qquad a^{-1}\in L^1(0,T^*).
\]
Let the initial Cauchy data agree with the FLRW data on the interior
ball $\U=\{0\le s<s_0\}$, with area radius
$r_0=S_\kappa(s_0)$. If
\begin{equation}\label{eq:causal-collapse-area}
 r_0>\Rcoll:=S_\kappa(\eta_*),
 \qquad \eta_*:=\int_0^{T^*}\frac{\dd t}{a(t)},
\end{equation}
then the central comoving geodesic $\gamma$, given by $s=0$, remains in
$\operatorname{int}D^+(\U)$ for $0<t<T^*$ and satisfies
\[
 \operatorname{Scal}_g(\gamma(t))\longrightarrow-\infty
 \qquad(t\to T^*).
\]
\end{corollary}

\begin{proof}
By Proposition~\ref{prop:causal-collapse}, it remains to show that $a(t)\to0$ as
$t\to T^*$. Since $\dot a<0$ and $a>0$ on $[0,T^*)$, the limit
$a_*:=\lim_{t\to T^*}a(t)$ exists. Suppose that $a_*>0$. By
\eqref{eq:contracting-energy},
\[
 \frac{\dd}{\dd t}(a^6E)
 =6a^5\dot a E+a^6\dot E
 =6\mu^2a^5\dot a\phi^2\le0.
\]
Recalling that $a(0)=1$, we obtain
\[
 E(t)\le\frac{E(0)}{a(t)^6}\le\frac{E(0)}{a_*^6}.
\]
Thus $E$, $\phi$, and $\dot\phi$ are bounded on $[0,T^*)$, while
$0<a_*\le a(t)\le1$. The Friedmann equation \eqref{eq:flrw-constraint-intro} also bounds $H$.
Moreover, \eqref{eq:flrw-constraint-intro} and
\eqref{eq:flrw-evolution-intro} give
\[
 \dot H=\frac{\ddot a}{a}-H^2=\frac{\kappa}{a^2}-\dot\phi^2.
\]
Together with $\dot a=aH$ and \eqref{eq:flrw-kg-intro}, this bounds the derivatives
of $a$, $H$, $\phi$, and $\dot\phi$, so these variables have finite
limits at $T^*$. Since $a_*>0$, local existence and uniqueness for the
ODE system extend the solution smoothly beyond $T^*$.
Moreover, $\kappa\le0$ and $E(t)\ge E(0)>0$ give
\[
 H(t)=-\sqrt{\frac{E(t)}3-\frac{\kappa}{a(t)^2}}
 \le-\sqrt{\frac{E(0)}3}<0.
\]
Since the limiting values of $a$ and $H$ are strictly positive and
strictly negative, respectively, the extended solution still satisfies
$a>0$ and $H<0$ for a short time after $T^*$. This contradicts the
definition of $T^*$ as the endpoint of the maximal regular contracting
interval. Hence $a(t)\to0$, and Proposition~\ref{prop:causal-collapse} applies
because $r_0>\Rcoll$ is equivalent to $s_0>\eta_*$.
\end{proof}

The results above give a sufficient condition $r_0>\Rtrap$ for
trapped surface formation under the assumptions of
Proposition~\ref{prop:causal-criterion}, and a sufficient condition $r_0>\Rcoll$
for curvature singularity formation under the assumptions of
Corollary~\ref{cor:causal-collapse}. We also require $r_0<\Rglue$ to extend the
data on the interior FLRW ball to complete untrapped initial data, as stated in Proposition~\ref{prop:completion}.
In the next section, we compare the three radius thresholds $\Rtrap$, $\Rcoll$, and $\Rglue$ and
prove the main results stated in \cref{sec:main-results}.

\section{Comparison of the radius thresholds}\label{sec:comparison}

The results of the previous two sections reduce the construction of
the initial data to a comparison of the radius thresholds. By
Proposition~\ref{prop:completion}, Proposition~\ref{prop:causal-criterion}, and Lemma~\ref{lem:touching}, an initial
FLRW ball with area radius $r_0$ satisfying
\begin{equation}\label{eq:trapping-gluing-window}
 \Rtrap<r_0<\Rglue
\end{equation}
can be extended to complete, initially untrapped, asymptotically flat
initial data whose future development contains a trapped surface.
Under the additional assumption \eqref{eq:kinetic-cone-intro} on the initial scalar data,
\begin{equation}\label{eq:comparison-kinetic-assumption}
 \Pi_0^2>2\mu^2\Phi_0^2,
\end{equation}
we prove that $\Rtrap<\Rglue$ and hence obtain Theorem~\ref{thm:general}.

\begin{proof}[Proof of Theorem~\ref{thm:general}]
Recall from \eqref{eq:horizon-coordinate} and
\eqref{eq:trapping-coordinates} that
\[
 \begin{gathered}
 s_{\mathrm{AH}}(t)
 =S_\kappa^{-1}\left(\frac1{\sqrt{\dot a(t)^2+\kappa}}\right),\\
 \widetilde{s}(t)=\eta(t)+s_{\mathrm{AH}}(t),
 \qquad \eta(t)=\int_0^t\frac{\dd\tau}{a(\tau)}.
 \end{gathered}
\]
With $E$ defined by \eqref{eq:scalar-energy}, the Friedmann equations
\eqref{eq:flrw-constraint-intro}--\eqref{eq:flrw-evolution-intro} give
\[
 \dot a^2+\kappa=\frac{a^2E}{3},
 \qquad
 \ddot a=\frac a3\bigl(\mu^2\phi^2-2\dot\phi^2\bigr).
\]
On the increasing branch of $S_\kappa$,
\eqref{eq:curvature-identities} and $\dot a<0$ imply
\[
 S_\kappa'(s_{\mathrm{AH}})
 =\sqrt{1-\kappa S_\kappa(s_{\mathrm{AH}})^2}
 =\frac{-\dot a}{\sqrt{\dot a^2+\kappa}}>0.
\]
Direct differentiation with respect to $t$ gives
\begin{equation}\label{eq:F-derivative}
 \begin{aligned}
 \dot{\widetilde{s}}(t)
 &=\dot\eta(t)+\dot s_{\mathrm{AH}}(t)\\
 &=\frac1{a(t)}+
   \frac{\frac{\dd}{\dd t}S_\kappa(s_{\mathrm{AH}}(t))}
        {S_\kappa'(s_{\mathrm{AH}}(t))}\\
 &=\frac1{a(t)}+
   \frac{-\dot a(t)\ddot a(t)}{(\dot a(t)^2+\kappa)^{3/2}}
   \frac{\sqrt{\dot a(t)^2+\kappa}}{-\dot a(t)}\\
 &=\frac1{a(t)}+\frac{\ddot a(t)}{\dot a(t)^2+\kappa}
 =\frac{a(t)E(t)+3\ddot a(t)}{a(t)^2E(t)}\\
 &=\frac{2\mu^2\phi(t)^2-\dot\phi(t)^2}{a(t)E(t)}.
 \end{aligned}
\end{equation}
This identity holds on any regular contracting interval.
By \eqref{eq:comparison-kinetic-assumption} and $a(0)=1$,
\begin{equation}\label{eq:initial-threshold-derivative}
 \dot{\widetilde{s}}(0)
 =\frac{2\mu^2\Phi_0^2-\Pi_0^2}{E_0}<0,
 \qquad E_0=E(0).
\end{equation}
Thus, for sufficiently small positive $t$,
$\widetilde{s}(t)<\widetilde{s}(0)$.
Since $\widetilde{s}(0)=s_{\mathrm{AH}}(0)$, the definition
$s_{\mathrm{trap}}=\inf_{0<t<T^*}\widetilde{s}(t)$ gives
\begin{equation}\label{eq:initial-trapping-comparison}
 s_{\mathrm{trap}}<s_{\mathrm{AH}}(0).
\end{equation}

The same assumption gives $E_0>3\mu^2\Phi_0^2$.
The definitions of $\Rglue$, $r_{\mathrm{AH}}$, and
$s_{\mathrm{AH}}$ in \eqref{eq:completion-condition},
\eqref{eq:apparent-area-radius}, and \eqref{eq:horizon-coordinate},
together with $a(0)=1$, therefore yield
\begin{equation}\label{eq:initial-gluing-horizon}
 \Rglue=\sqrt{\frac3{E_0}}
 =r_{\mathrm{AH}}(0)=S_\kappa(s_{\mathrm{AH}}(0)).
\end{equation}
For $\kappa\le0$, $S_\kappa$ is strictly increasing on $[0,\infty)$.
For $\kappa>0$, the condition $\kappa<E_0/3$ gives
\[
 s_{\mathrm{AH}}(0)
 =S_\kappa^{-1}\left(\sqrt{\frac3{E_0}}\right)
 \in\left(0,\frac{\pi}{2\sqrt\kappa}\right),
\]
so $S_\kappa$ is also strictly increasing on
$[0,s_{\mathrm{AH}}(0)]$.
Since $\Rtrap=S_\kappa(s_{\mathrm{trap}})$,
\eqref{eq:initial-trapping-comparison} and
\eqref{eq:initial-gluing-horizon} imply $\Rtrap<\Rglue$.

Now choose any $r_0$ with $\Rtrap<r_0<\Rglue$.
By Proposition~\ref{prop:completion} and Lemma~\ref{lem:touching}, the interior FLRW data in
\eqref{eq:completed-data-intro} for $0\le r<r_0$ extend to smooth,
complete, asymptotically flat initial data on $\R^3$ with no closed
weakly outer trapped surfaces. Since $r_0>\Rtrap$,
Proposition~\ref{prop:causal-criterion} gives a trapped sphere in $D^+(\U)$.
Finally, for every null vector $L$, the Einstein equation
\eqref{eq:ekg-einstein} gives
\[
 \Ric(L,L)=2(L\phi)^2\ge0.
\]
The maximal Cauchy development has a noncompact Cauchy hypersurface.
The Penrose singularity theorem \cite{Penrose1965} therefore implies
future null geodesic incompleteness.
\end{proof}

To prove Theorem~\ref{thm:small-radius}, we consider sufficiently negative
curvature. Keep $\mu>0$, $\Phi_0\in\R$, and $\Pi_0\ne0$ fixed, and write
\begin{equation}\label{eq:negative-family-data}
 \kappa=-A^2,\qquad
 a(0)=1,\qquad \phi(0)=\Phi_0,\qquad \dot\phi(0)=\Pi_0,
\end{equation}
where $A>0$. Thus $\kappa\to-\infty$ means $A\to\infty$.

Since $\dot a(t)<0$ on the maximal regular contracting interval
$[0,T^*)$, the Friedmann equation \eqref{eq:flrw-constraint-intro} gives
\begin{equation}\label{eq:negative-friedmann}
 \dot a(t)=-\sqrt{A^2+\frac{a(t)^2E(t)}3}\le-A.
\end{equation}
Recalling that $a(0)=1$, we have $0<a(t)\le1$ and
$a(t)\le1-At$ on $[0,T^*)$. Hence
\begin{equation}\label{eq:negative-collapse-time}
 0<T^*\le\frac1A.
\end{equation}
The proof of Corollary~\ref{cor:causal-collapse} shows that, for $\kappa\le0$,
a finite endpoint of the maximal regular contracting interval implies
$a(t)\to0$ as $t\to T^*$. This part of that proof does not require
finite conformal time and applies here because $\kappa<0$.

To control the scalar field $\phi(t)$, set
\begin{equation}\label{eq:p-def}
 p(t):=a(t)^3\dot\phi(t).
\end{equation}
The Klein--Gordon equation \eqref{eq:flrw-kg-intro} then becomes
\begin{equation}\label{eq:p-evolution}
 \dot p(t)=-\mu^2a(t)^3\phi(t),\qquad p(0)=\Pi_0.
\end{equation}
We next obtain positive lower and upper bounds for $|p(t)|$
throughout $[0,T^*)$.

\begin{lemma}\label{lem:momentum-control}
Let $A$ and $p$ be defined by \eqref{eq:negative-family-data} and
\eqref{eq:p-def}. For all sufficiently large $A>0$, $p$ has the same
sign as $\Pi_0$ and satisfies
\begin{equation}\label{eq:p-bounds}
 \frac{|\Pi_0|}{2}\le |p(t)|\le2|\Pi_0|
 \qquad(0\le t<T^*).
\end{equation}
\end{lemma}

\begin{proof}
Since $p(0)=\Pi_0$, we have $|p(t)|<2|\Pi_0|$ for sufficiently
small $t\ge0$. Let $[0,t_*)$ be the maximal subinterval of $[0,T^*)$
on which $|p(t)|\le2|\Pi_0|$.
For $0\le t<t_*$, \eqref{eq:p-def} and
\eqref{eq:negative-friedmann} give
\[
 \begin{aligned}
 |\phi(t)|
 &\le |\Phi_0|+\int_0^t\frac{|p(\tau)|}{a(\tau)^3}\,\dd\tau\\
 &\le |\Phi_0|+\frac{2|\Pi_0|}{A}
       \int_{a(t)}^1\frac{\dd x}{x^3}\\
 &=|\Phi_0|+\frac{|\Pi_0|}{A}\bigl(a(t)^{-2}-1\bigr).
 \end{aligned}
\]
Here we changed variables from $\tau$ to $x=a(\tau)$ and used
$1/|\dot a|\le1/A$. Integrating \eqref{eq:p-evolution} and changing
variables in the same way, we obtain
\[
 \begin{aligned}
 |p(t)-\Pi_0|
 &\le\mu^2\int_0^t a(\tau)^3|\phi(\tau)|\,\dd\tau\\
 &\le\frac{\mu^2}{A}\int_{a(t)}^1
 x^3\left[|\Phi_0|+\frac{|\Pi_0|}{A}(x^{-2}-1)\right]\dd x\\
 &\le\frac{\mu^2|\Phi_0|}{4A}
       +\frac{\mu^2|\Pi_0|}{4A^2}.
 \end{aligned}
\]
For fixed $\mu,\Phi_0,\Pi_0$, choosing $A$ sufficiently large gives
\begin{equation}\label{eq:p-loss}
 |p(t)-\Pi_0|\le\frac{|\Pi_0|}{2}
 \qquad(0\le t<t_*).
\end{equation}
If $t_*<T^*$, continuity gives $|p(t_*)|\le3|\Pi_0|/2$.
Thus $|p(t)|<2|\Pi_0|$ also holds for a short time after $t_*$,
contradicting the maximality of $t_*$. Hence $t_*=T^*$.
It follows from \eqref{eq:p-loss} that $p$ has the same sign as
$\Pi_0$ and
\[
 \frac{|\Pi_0|}{2}\le |p(t)|\le\frac{3|\Pi_0|}{2}
 \qquad(0\le t<T^*),
\]
which proves \eqref{eq:p-bounds}.
\end{proof}

Next we show that $a^{-1}$ is integrable on $(0,T^*)$ for sufficiently
large $A$. This is an important hypothesis for curvature blow-up in
Proposition~\ref{prop:causal-collapse}. We also estimate $\Rcoll$ and show that
$\Rcoll\to0$ as $A\to\infty$.

\begin{lemma}\label{lem:collapse-distance}
For all sufficiently large $A>0$, $a^{-1}$ is integrable on $(0,T^*)$,
so $\eta_*:=\int_0^{T^*}a(t)^{-1}\dd t$ is finite, and
\begin{equation}\label{eq:collapse-radius-bound}
 \Rcoll:=S_\kappa(\eta_*)=\frac1A\sinh(A\eta_*)
 \le\frac{3^{1/4}}{\sqrt{A|\Pi_0|}}\longrightarrow0
 \qquad(A\to\infty).
\end{equation}
\end{lemma}

\begin{proof}
The Friedmann equation \eqref{eq:flrw-constraint-intro} gives
\begin{equation}\label{eq:negative-friedmann-squared}
 \dot a^2=A^2+\frac{a^2E}{3}.
\end{equation}
Since $E=\dot\phi^2+\mu^2\phi^2\ge\dot\phi^2$ and $p=a^3\dot\phi$,
we have
\[
 E\ge\frac{p^2}{a^6}.
\]
For sufficiently large $A$, Lemma~\ref{lem:momentum-control} gives
$|p|\ge|\Pi_0|/2$, and hence
\[
 E\ge\frac{\Pi_0^2}{4a^6}.
\]
Substituting this into \eqref{eq:negative-friedmann-squared}, we obtain
\begin{equation}\label{eq:collapse-speed-bound}
 \dot a^2\ge A^2+\frac{\Pi_0^2}{12a^4}.
\end{equation}
As $a(t)$ decreases from $1$ to $0$ on $[0,T^*)$, we change variables
from $t$ to $a$, apply \eqref{eq:collapse-speed-bound}, and then
set $x=a^2$ to obtain
\begin{equation}\label{eq:eta-star-bound}
 \begin{aligned}
 \eta_*
 &=\int_0^{T^*}\frac{\dd t}{a(t)}
 =\int_0^1\frac{\dd a}{a|\dot a|}\\
 &\le\frac1A\int_0^1
       \frac{a\,\dd a}{\sqrt{a^4+\Pi_0^2/(12A^2)}}\\
 &=\frac1{2A}\int_0^1
       \frac{\dd x}{\sqrt{x^2+\Pi_0^2/(12A^2)}}\\
 &=\frac1{2A}\operatorname{arsinh}
       \left(\frac{2\sqrt3 A}{|\Pi_0|}\right)<\infty.
 \end{aligned}
\end{equation}
Thus $a^{-1}$ is integrable on $(0,T^*)$, and $\Rcoll$ is well defined.
Since $S_{-A^2}(s)=A^{-1}\sinh(As)$ by \eqref{eq:S-kappa-intro},
the estimate \eqref{eq:eta-star-bound} gives
\begin{equation}\label{eq:collapse-radius-intermediate}
 \Rcoll=\frac1A\sinh(A\eta_*)
 \le\frac1A\sinh\left(\frac12\operatorname{arsinh}
                       \frac{2\sqrt3 A}{|\Pi_0|}\right).
\end{equation}
For $x\ge0$,
\[
 \sinh\left(\frac12\operatorname{arsinh}x\right)
 =\sqrt{\frac{\sqrt{1+x^2}-1}{2}}\le\sqrt{x/2}.
\]
Applying this inequality to \eqref{eq:collapse-radius-intermediate}
gives
\[
 \Rcoll\le\frac1A\sqrt{\frac{\sqrt3 A}{|\Pi_0|}}
 =\frac{3^{1/4}}{\sqrt{A|\Pi_0|}}\longrightarrow0
 \qquad(A\to\infty),
\]
which proves \eqref{eq:collapse-radius-bound}.
\end{proof}

Next we show that $\Rtrap$ equals $\Rcoll$ when $A$ is sufficiently large.

\begin{lemma}\label{lem:threshold-equality}
For all sufficiently large $A>0$,
\begin{equation}\label{eq:trap-equals-collapse}
 \Rtrap=\Rcoll.
\end{equation}
\end{lemma}

\begin{proof}
Since $\Rtrap=S_\kappa(s_{\mathrm{trap}})$ and
$\Rcoll=S_\kappa(\eta_*)$, it suffices to prove
\begin{equation}\label{eq:geodesic-trap-equals-collapse}
 s_{\mathrm{trap}}=\eta_*.
\end{equation}
Recall from \eqref{eq:trapping-coordinates} that
\[
 \begin{gathered}
 s_{\mathrm{trap}}=\inf_{0<t<T^*}\widetilde{s}(t),
 \qquad \widetilde{s}(t)=\eta(t)+s_{\mathrm{AH}}(t),\\
 \eta(t)=\int_0^t\frac{\dd\tau}{a(\tau)}.
 \end{gathered}
\]
We first consider the values of $\widetilde{s}$ at the two endpoints.
At $t=0$,
\begin{equation}\label{eq:threshold-initial-value}
 \widetilde{s}(0)=\eta(0)+s_{\mathrm{AH}}(0)=s_{\mathrm{AH}}(0).
\end{equation}
Recall that $S_\kappa^{-1}$ is increasing for $\kappa=-A^2<0$.
By \eqref{eq:horizon-coordinate} and \eqref{eq:flrw-constraint-intro},
\[
 \begin{aligned}
 s_{\mathrm{AH}}(t)
 &=S_\kappa^{-1}\left(\frac1{\sqrt{\dot a(t)^2+\kappa}}\right)\\
 &=S_\kappa^{-1}\left(
     \frac{\sqrt3}{a(t)\sqrt{\dot\phi(t)^2+\mu^2\phi(t)^2}}\right)\\
 &\le S_\kappa^{-1}\left(\frac{\sqrt3}{a(t)|\dot\phi(t)|}\right)
 =S_\kappa^{-1}\left(\frac{\sqrt3\,a(t)^2}{|p(t)|}\right).
 \end{aligned}
\]
By Lemma~\ref{lem:momentum-control}, $|p(t)|\ge|\Pi_0|/2$, so
\[
 0<s_{\mathrm{AH}}(t)
 \le S_\kappa^{-1}\left(\frac{2\sqrt3\,a(t)^2}{|\Pi_0|}\right).
\]
Since $a(t)\to0$ and $S_\kappa^{-1}$ is continuous at zero, we obtain
\begin{equation}\label{eq:horizon-coordinate-limit}
 s_{\mathrm{AH}}(t)\longrightarrow0
 \qquad(t\to T^*).
\end{equation}
Together with Lemma~\ref{lem:collapse-distance}, this gives
\begin{equation}\label{eq:threshold-terminal-limit}
 \lim_{t\to T^*}\widetilde{s}(t)
 =\lim_{t\to T^*}\eta(t)+\lim_{t\to T^*}s_{\mathrm{AH}}(t)
 =\eta_*.
\end{equation}
Thus $\widetilde{s}$ extends continuously to $[0,T^*]$, with
\begin{equation}\label{eq:threshold-endpoint-values}
 \widetilde{s}(0)=s_{\mathrm{AH}}(0),
 \qquad \widetilde{s}(T^*)=\eta_*.
\end{equation}

Suppose that $\widetilde{s}$ attains its minimum at some
$t\in(0,T^*)$. Then $\dot{\widetilde{s}}(t)=0$ and
$\ddot{\widetilde{s}}(t)\ge0$. By \eqref{eq:F-derivative},
\[
 \dot{\widetilde{s}}(t)
 =\frac{2\mu^2\phi(t)^2-\dot\phi(t)^2}{E(t)a(t)}=0,
\]
so $2\mu^2\phi(t)^2-\dot\phi(t)^2=0$.
When we differentiate this quotient at the minimum point, the term
containing the derivative of the denominator vanishes. Therefore,
\[
 \begin{aligned}
 \ddot{\widetilde{s}}(t)
 &=\frac{4\mu^2\phi(t)\dot\phi(t)
          -2\dot\phi(t)\ddot\phi(t)}{E(t)a(t)}\\
 &=\frac{6\mu^2\phi(t)\dot\phi(t)
          +6H(t)\dot\phi(t)^2}{E(t)a(t)},
 \end{aligned}
\]
where the last equality uses the Klein--Gordon equation
\eqref{eq:flrw-kg-intro}.
At this point, $\dot\phi(t)^2=2\mu^2\phi(t)^2$, and hence
$|\phi(t)\dot\phi(t)|=\dot\phi(t)^2/(\sqrt2\mu)$.
The Friedmann equation \eqref{eq:flrw-constraint-intro} also gives
\[
 H(t)=-\sqrt{\frac{E(t)}3-\frac{\kappa}{a(t)^2}}
 =-\sqrt{\frac{E(t)}3+\frac{A^2}{a(t)^2}}
 \le-\frac{A}{a(t)}.
\]
Consequently, if $A\ge\mu>0$, then $a(t)\le1$ gives
\[
 \begin{aligned}
 \ddot{\widetilde{s}}(t)
 &\le\left(3\sqrt2\mu-\frac{6A}{a(t)}\right)
       \frac{\dot\phi(t)^2}{E(t)a(t)}\\
 &\le(3\sqrt2\mu-6A)\frac{\dot\phi(t)^2}{E(t)a(t)}<0.
 \end{aligned}
\]
The last inequality is strict because
$\dot\phi(t)\ne0$ by \eqref{eq:p-bounds}.
This contradicts $\ddot{\widetilde{s}}(t)\ge0$.
Thus the minimum of $\widetilde{s}$ is attained at $t=0$ or $t=T^*$.

Finally,
\[
 S_\kappa(s_{\mathrm{AH}}(0))=\sqrt{\frac3{E_0}}>0,
 \qquad
 S_\kappa(\eta_*)=\Rcoll\longrightarrow0
 \quad(A\to\infty).
\]
The first radius is independent of $A$. Since $S_\kappa$ is strictly
increasing, for sufficiently large $A$ we have
\[
 \widetilde{s}(0)=s_{\mathrm{AH}}(0)>\eta_*=\widetilde{s}(T^*).
\]
The endpoint values are approached from $(0,T^*)$, so
\[
 s_{\mathrm{trap}}=\inf_{0<t<T^*}\widetilde{s}(t)
 =\widetilde{s}(T^*)=\eta_*.
\]
This proves \eqref{eq:geodesic-trap-equals-collapse}, and hence
$\Rtrap=S_\kappa(s_{\mathrm{trap}})=S_\kappa(\eta_*)=\Rcoll$,
as claimed in \eqref{eq:trap-equals-collapse}.
\end{proof}

Combining Lemmas~\ref{lem:momentum-control}, \ref{lem:collapse-distance}, and~\ref{lem:threshold-equality},
we can now prove Theorem~\ref{thm:small-radius}.

\begin{proof}[Proof of Theorem~\ref{thm:small-radius}]
Set $\kappa=-A^2$ and choose $A>0$ sufficiently large that
Lemmas~\ref{lem:momentum-control}, \ref{lem:collapse-distance}, and~\ref{lem:threshold-equality}
apply. These lemmas give
\[
 \Rcoll=\Rtrap\longrightarrow0
 \qquad(A\to\infty),
\]
whereas the gluing threshold
\[
 \Rglue=\min\left\{\sqrt{\frac3{E_0}},\frac1{\mu|\Phi_0|}\right\}>0
\]
depends only on the fixed data $\mu,\Phi_0,\Pi_0$.
Therefore, for sufficiently large $A$,
\[
 \Rcoll=\Rtrap<\Rglue.
\]
Choose any area radius $r_0$ with
\[
 \Rtrap=\Rcoll<r_0<\Rglue.
\]
By Proposition~\ref{prop:completion} and Lemma~\ref{lem:touching}, the interior FLRW data in
\eqref{eq:completed-data-intro} on $\U=\{0\le r<r_0\}$ extend to
smooth, complete, initially untrapped, spherically symmetric,
asymptotically flat data on $\R^3$. The extension agrees with
Schwarzschild data outside a larger ball of area radius $R_0>r_0$.
Since $r_0>\Rtrap$, Proposition~\ref{prop:causal-criterion} gives a trapped sphere
in $D^+(\U)$.

Moreover, \eqref{eq:negative-collapse-time} and
Lemma~\ref{lem:collapse-distance} give $T^*<\infty$ and
$a^{-1}\in L^1(0,T^*)$. Since $\kappa<0$ and $r_0>\Rcoll$,
Corollary~\ref{cor:causal-collapse} gives
\[
 \operatorname{Scal}_g(\gamma(t))\longrightarrow-\infty
 \qquad(t\to T^*)
\]
along the central comoving geodesic $\gamma$.
Its proper time is $t$, so it has finite future proper length.
If $\gamma$ extended through $T^*$ in the smooth maximal Cauchy
development, scalar curvature would be finite and continuous at its
endpoint, contrary to this limit. Thus the Cauchy development is
future timelike geodesically incomplete.
\end{proof}

The comparison of the radius thresholds also gives initial data with
arbitrarily small compactness.

\begin{proof}[Proof of Corollary~\ref{cor:small-initial-compactness}]
Fix $\mu>0$, $r_0>0$, and $\delta\in(0,1/2)$.
Take $R_0=2r_0$, set $\Phi_0=0$, and choose $\Pi_0$ such that
\[
 0<|\Pi_0|<\frac{\sqrt{3\delta}}{2r_0}.
\]
Then
\[
 r_0<R_0<\frac{\sqrt3}{|\Pi_0|}=\Rglue,
 \qquad
 \left(\frac{R_0}{\Rglue}\right)^2
 =\frac{\Pi_0^2R_0^2}{3}<\delta.
\]
Keep $\Phi_0$ and $\Pi_0$ fixed and take $\kappa=-A^2$ sufficiently
negative. By Lemmas~\ref{lem:collapse-distance} and~\ref{lem:threshold-equality},
we then have
\[
 \Rtrap=\Rcoll<r_0.
\]
Now Proposition~\ref{prop:small-initial-mass-ratio} gives smooth, complete,
initially untrapped, spherically symmetric, asymptotically flat
initial data satisfying \eqref{eq:general-initial-mass-ratio}, so
\[
 \sup_\Sigma\frac{2\HM}{r}
 \le\left(\frac{R_0}{\Rglue}\right)^2
 =\frac{\Pi_0^2R_0^2}{3}<\delta.
\]
These data agree with the same FLRW data on the interior ball.
Since $r_0>\Rtrap=\Rcoll$, Proposition~\ref{prop:causal-criterion} and Corollary~\ref{cor:causal-collapse}
give trapped surface formation and curvature blow-up in the future
development.
\end{proof}

\begin{figure}[!htbp]
\centering
\begin{tikzpicture}[x=1.08cm,y=1.08cm,>=Latex,
  every node/.style={font=\small},line cap=round,line join=round]
  \coordinate (O) at (0,0);
  \coordinate (B) at (4,0);
  \coordinate (I) at (7.6,0);
  \coordinate (P) at (0,-1.62);
  \coordinate (V) at (0,-4);
  \coordinate (S) at (0,2.25);
  \coordinate (Q) at (1.75,2.25);
  \fill[black!3] (O) -- (I) -- (5.35,2.25)
    .. controls (4.1,2.65) and (2.75,2.25) .. (Q) -- (S) -- cycle;
  \fill[black!7] (V) -- (B) -- (Q) -- (S) -- cycle;
  \draw[black!55,densely dotted,line width=1.05pt] (Q)
    .. controls (2.75,2.25) and (4.1,2.65) .. (5.35,2.25) -- (I);
  \draw[thick] (V) -- (S);
  \draw[dashed,thick] (B) -- (Q);
  \draw[dashed,thick] (V) -- (B);
  \draw[black!45,line width=1pt] (O) -- (3.69,0);
  \draw[line width=1.55pt] (P) -- (.22,-1.62)
    .. controls (1.37,-1.62) and (2.54,0) .. (3.69,0) -- (I);
  \draw[line width=1.2pt,decorate,
    decoration={zigzag,segment length=4pt,amplitude=2.3pt}]
    (0.02,2.25) -- (1.73,2.25);
  \fill (O) circle (1.35pt);
  \fill (P) circle (1.6pt);
  \fill (B) circle (1.5pt);
  \fill (I) circle (1.5pt);
  \node[anchor=south west,font=\footnotesize] at (-.05,2.33) {curvature blow-up};
  \draw[decorate,decoration={brace,amplitude=4pt,raise=3pt}]
    (O) -- (S) node[midway,left=10pt,font=\footnotesize] {$T^*$};
  \node at (1.25,1.0) {$D^+(\U)$};
  \node at (.73,-2.22) {$D^-(\U)$};
  \draw[thin,black!65,decorate,decoration={brace,mirror,amplitude=3pt,raise=1.5pt}]
    (O) -- (B) node[midway,below=5pt,text=black] {$\U$};
  \node[anchor=south] at (1,.00) {$\Sigma$};
  \node[anchor=south west,font=\footnotesize] at (3.7,.03) {$r=r_0$};
  \node[below,font=\footnotesize] at (5.8,-.05)
    {common exterior of $\Sigma$ and $\widetilde\Sigma$};
  \node[right] at (I) {$i^0$};
  \node[anchor=north west] at (0.68,-.89) {$\widetilde\Sigma$};
  \draw[decorate,decoration={brace,amplitude=4pt,raise=3pt}]
    (P) -- (O) node[midway,left=10pt,font=\footnotesize] {$r_0$};
\end{tikzpicture}
\caption{The deformed slice $\widetilde\Sigma$ lies to the past
of the original slice $\Sigma$ near the center. The two
slices coincide near $\partial\U$ and throughout the exterior. The interior of $\widetilde\Sigma$ lies in the past domain of dependence $D^-(\mathcal{U})$ of $\mathcal{U}$.}
\label{fig:past-slice}
\end{figure}
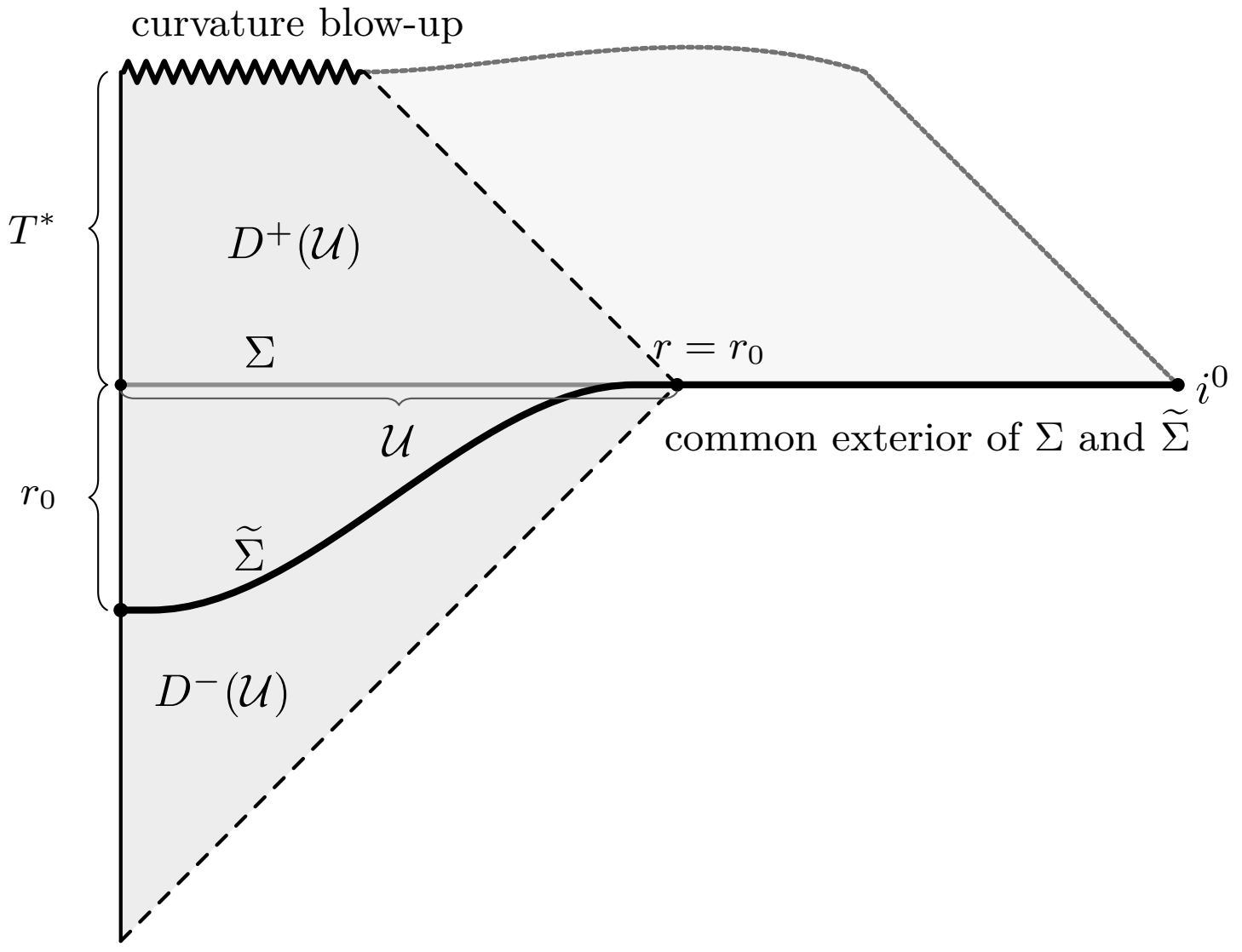

In the construction above, let $T_{\mathrm{center}}(\Sigma)$ denote the
proper time from the initial hypersurface $\Sigma$ to curvature blow-up
along the central timelike geodesic. Since $t$ is proper time on this
geodesic, \eqref{eq:negative-collapse-time} gives
\begin{equation}\label{eq:original-central-time}
 T_{\mathrm{center}}(\Sigma)=T^*\le A^{-1}.
\end{equation}
For fixed $\mu>0$ and $r_0>0$, letting $\kappa=-A^2\to-\infty$ gives
\begin{equation}\label{eq:original-central-mass-scale}
 \mu T_{\mathrm{center}}(\Sigma)\longrightarrow0.
\end{equation}
Thus the fixed scalar-field mass is small in units set by \(T_{\mathrm{center}}(\Sigma)\).
If we take $T_{\mathrm{center}}(\Sigma)$ as the characteristic time of the
spacetime, \eqref{eq:original-central-mass-scale} can be interpreted as the limit $\mu\to0$ for the rescaled Klein--Gordon equation. Our aim is
to allow any prescribed positive value of
the dimensionless mass $\mu T_{\mathrm{center}}$.
Moving the central part of $\Sigma$ to the past
gives the proper time $T_{\mathrm{center}}$ a lower bound independent of $\kappa$.
The new Cauchy hypersurface $\tilde{\Sigma}$ will coincide with $\Sigma$ outside $\U$.
Its deformed portion lies in the past domain of dependence $D^-(\U)$, as shown in \cref{fig:past-slice}.
We first control the compactness throughout $D^-(\U)$.

\begin{lemma}\label{lem:past-compactness}
For the data constructed in the proof of
Corollary~\ref{cor:small-initial-compactness},
\[
 \sup_{D^-(\U)}\frac{2\HM}{r}<\delta .
\]
\end{lemma}

\begin{proof}
Since $\kappa=-A^2<0$, the Friedmann equation
\eqref{eq:flrw-constraint-intro} and $H(0)<0$ imply $H(t)<0$ throughout
the regular past of the FLRW solution. By \eqref{eq:scalar-energy}
and \eqref{eq:contracting-energy},
\[
 \begin{aligned}
  E(t)      &= \dot\phi(t)^2+\mu^2\phi(t)^2,\\[2pt]
  \dot E(t) &= -6H(t)\dot\phi(t)^2 \ge 0.
 \end{aligned}
\]
Since $\Phi_0=0$ in this construction, it follows that
\[
 E(t)\le E(0)=\Pi_0^2 \qquad (t\le 0).
\]

From each sphere of symmetry in the interior of $D^-(\U)$, the
future outgoing radial null geodesic reaches $\U$ on $\Sigma$. Let $\lambda$ be an affine parameter along such a geodesic,
and set
\[
 \ell:=\frac{\dd}{\dd\lambda},
 \qquad
 \theta:=\frac{2\ell(r)}{r},
\]
where $\theta$ is the null expansion.  By spherical symmetry and \eqref{eq:ekg-einstein},
the Raychaudhuri equation reads
\[
 \begin{aligned}
  \ell(\theta) &= -\frac12\theta^2-\Ric(\ell,\ell)\\[2pt]
               &= -\frac12\theta^2-2(\ell\phi)^2 \le 0.
 \end{aligned}
\]
The outgoing expansion is positive on $\U$, so it is positive along
the preceding segment of the ray. The area radius $r$ in
\eqref{eq:area-radius} therefore increases towards $\U$, and $r\le r_0$
throughout $D^-(\U)$. Applying the FLRW mass formula
\eqref{eq:ball-masses} on each constant-$t$ slice yields
\[
 \begin{aligned}
  \frac{2\HM}{r} = \frac{E(t)r^2}{3}\le \frac{\Pi_0^2 r_0^2}{3}<\delta ,
 \end{aligned}
\]
where the final bound follows from the choice of $\Pi_0$ in the proof
of Corollary~\ref{cor:small-initial-compactness}.
\end{proof}

\begin{proposition}\label{prop:central-proper-time}
Fix $\mu>0$, $r_0>0$ and $\delta\in(0,1/2)$. Consider the data on
$\Sigma$ constructed in the proof of Corollary~\ref{cor:small-initial-compactness}
with $\kappa=-A^2$ and $A$ sufficiently large. Their maximal Cauchy
development admits a smooth, complete, spherically symmetric,
asymptotically flat Cauchy hypersurface $\widetilde\Sigma$ that agrees
with $\Sigma$ outside $\U$ and satisfies
\begin{equation}\label{eq:shifted-compactness}
 \sup_{\widetilde\Sigma}\frac{2\HM}{r}<\delta .
\end{equation}
This hypersurface is initially untrapped and has a trapped sphere in
its future. The proper time from $\widetilde\Sigma$ to curvature
blow-up along the central timelike geodesic satisfies
\begin{equation}\label{eq:shifted-central-time}
 T_{\mathrm{center}}(\widetilde\Sigma)=r_0+T^*
 \in(r_0,r_0+A^{-1}].
\end{equation}
\end{proposition}

\begin{proof}
By the Friedmann equation \eqref{eq:flrw-constraint-intro} and the
past energy bound from the proof of Lemma~\ref{lem:past-compactness},
\begin{equation}\label{eq:dota}
 A\le -\dot a(t)=\sqrt{A^2+\frac{a(t)^2E(t)}{3}}
   \le A+a(t)\sqrt{\frac{E(0)}{3}} .
\end{equation}
The upper bound in \eqref{eq:dota} ensures that the solution remains regular back to $t=-r_0$. Integrating the
lower bound in \eqref{eq:dota} and using $a(0)=1$, we obtain \[a(t)\ge 1+A|t|\] on $[-r_0,0]$.

Recall the definition of conformal time in \eqref{eq:trapping-coordinates}.
Then
\[
 -\eta(-r_0)
 =\int_{-r_0}^{0}\frac{\dd\tau}{a(\tau)}
 \le \int_{-r_0}^{0}\frac{\dd\tau}{1+A|\tau|}
 =\frac{\log(1+Ar_0)}{A}.
\]
By \eqref{eq:U-geodesic} and \eqref{eq:inverse-radius}, the comoving
radius of $\U$ is
\[
 s_0=S_\kappa^{-1}(r_0)
 =\frac{\operatorname{arsinh}(Ar_0)}{A}
 >\frac{\log(1+Ar_0)}{A}.
\]
Thus \[-\eta(-r_0)<s_0.\] Applying the radial null-ray argument of
\eqref{eq:Dplus-slice} towards the past, we see that the point
$(t,s)=(-r_0,0)$ on the central timelike geodesic lies in $D^-(\U)$.
The strict inequality allows us to join this point to
$\Sigma\setminus\U$ by a smooth spacelike hypersurface that agrees
with $\Sigma$ near the boundary of $\U$. The resulting Cauchy
hypersurface $\widetilde\Sigma$ lies in $D^-(\U)$ on its interior
portion and coincides with $\Sigma$ on the exterior.

By Lemma~\ref{lem:past-compactness}, the supremum of $2\HM/r$ on the
interior portion is less than $\delta$. The same bound holds on
$\Sigma\setminus\U$ by Corollary~\ref{cor:small-initial-compactness}. Hence
\[
 \sup_{\widetilde\Sigma}\frac{2\HM}{r}<\delta .
\]
The outgoing expansion is positive on $\tilde{\Sigma}$, so
Lemma~\ref{lem:touching} shows that $\tilde{\Sigma}$ is initially untrapped. The trapped sphere
in $D^+(\U)$ forms in the future of $\widetilde\Sigma$.

By construction, the proper time from the center of $\widetilde\Sigma$
to curvature blow-up along the central timelike geodesic is $r_0+T^*$.
Since $0<T^*\le A^{-1}$ by \eqref{eq:original-central-time}, we obtain
\[
 T_{\mathrm{center}}(\widetilde\Sigma)=r_0+T^*
 \in (r_0,r_0+A^{-1}]. \qedhere
\]
\end{proof}

With $\mu$ and $r_0$ fixed, \eqref{eq:shifted-central-time} gives
$\mu T_{\mathrm{center}}(\widetilde\Sigma)\to\mu r_0$ as $A\to\infty$.
The choice of $r_0$ in Corollary~\ref{cor:small-initial-compactness} allows
this limit to be any prescribed positive number.

\FloatBarrier

\section{Discussion}\label{sec:discussion}

The trapped surface formation established in
Theorems~\ref{thm:general} and~\ref{thm:small-radius} persists under sufficiently small
perturbations of the initial data. For this perturbation argument,
we introduce Cartesian coordinates $x$ on $\R^3$ associated with
the original spherical coordinates, so that the original radial
coordinate satisfies $r=|x|$. We express both the original and
perturbed data in these fixed coordinates and require the perturbed
data to be smooth and asymptotically flat and to satisfy the Einstein
constraint equations.
If the perturbed data are sufficiently close to our constructed data
in weighted Sobolev norms of sufficiently high order, the coordinate spheres
$\{|x|=\mathrm{constant}\}$ retain positive outward null expansion
at every point. A comparison argument as in the proof of
Lemma~\ref{lem:touching} then excludes closed weakly outer trapped surfaces
on the initial hypersurface. This comparison does not
require spherical symmetry. The perturbed data are also asymptotically flat.
By Cauchy stability \cite{HawkingEllis1973,HughesKatoMarsden1977,FriedrichRendall2000},
if the perturbed data are sufficiently close to our constructed data
in suitable weighted Sobolev norms, the corresponding Cauchy developments
remain close to the original Cauchy development on a region with compact
closure containing the trapped surfaces, after a suitable
identification of the spacetimes.
The metrics and their first derivatives are close enough that both
null expansions remain negative on the corresponding surface, so
trapped surfaces still form.

Thus each initial data set in Theorems~\ref{thm:general} and~\ref{thm:small-radius}
has an open neighborhood in the space of smooth asymptotically flat data satisfying the
constraint equations, with the topology defined by a suitable weighted Sobolev norm. The data in this neighborhood are
initially untrapped but give rise to trapped surface formation in the future development. These perturbations need not preserve spherical symmetry.

For the perturbed data above, a trapped surface forms in the future development,
so Penrose's theorem \cite{Penrose1965} gives future null geodesic incompleteness.
Now consider curvature blow-up for the perturbed data. For $\mu=0$,
Beyer-Oliynyk \cite{BeyerOliynyk2024} proved that scalar curvature blow-up
is stable under possibly nonsymmetric perturbations of spatially flat FLRW data,
provided that the perturbations are small in suitable Sobolev norms on a ball.
Their method and stability of the curvature blow-up may extend to the solutions constructed in
Theorem~\ref{thm:small-radius}, for which $\mu>0$ and $\kappa<0$.

The evolution outside $D^+(\U)$ is not fully described by our
analysis. Here $\U$ is the interior FLRW ball on the initial Cauchy
hypersurface. For the original data, the evolution inside $D^+(\U)$
coincides with the FLRW solution. The Schwarzschild exterior on the
initial Cauchy hypersurface determines the solution in its own
domain of dependence, but the entire exterior region outside $D^+(\U)$ is
not fully described. In particular, our analysis does not determine
the full event horizon of the solution.

Beyond the initial data constructed in this paper and their small
perturbations, we would like to find more general sufficient
conditions for trapped surface formation in the massive
Einstein--Klein--Gordon system that also apply to inhomogeneous data.
Such conditions should not require a smallness assumption on the
mass $\mu$ of the scalar field. A crucial issue in finding such
conditions is to control the mass term in the Klein--Gordon
equation, which requires estimates for the scalar field itself,
not only its derivatives.

\medskip
\noindent School of Mathematical Sciences, Peking University\\*
\textit{Email:} yxye.pku.edu.cn@stu.pku.edu.cn

\end{document}